\documentclass[11pt]{article} 

\usepackage[utf8]{inputenc} 

\usepackage{geometry} 
\usepackage{graphicx} 

\usepackage[parfill]{parskip} 

\usepackage{booktabs} 
\usepackage{array} 
\usepackage{paralist} 
\usepackage{verbatim} 
\usepackage{subfig} 

\usepackage{amsmath}
\usepackage{amssymb}
\usepackage{amsthm}
\usepackage{mathtools}
\usepackage{soul}

\begingroup
    \makeatletter
    \@for\theoremstyle:=definition,remark,plain\do{%
        \expandafter\g@addto@macro\csname th@\theoremstyle\endcsname{%
            \addtolength\thm@preskip\parskip
            }%
        }
\endgroup

\usepackage{algorithmic}
\usepackage{algorithm}

\usepackage{cite}
\usepackage{fullpage}

\usepackage{algorithmic}
\usepackage{algorithm}
\floatname{algorithm}{Method}

\usepackage{color}

\newtheorem{proposition}{Proposition}[section]
\newtheorem{assumption}{Assumption}[section]
\newtheorem{lemma}{Lemma}[section]
\newtheorem{corollary}{Corollary}[section]

\newtheorem{theorem}{Theorem}[section]
\theoremstyle{remark}

\newtheorem{remark}{Remark}[section]

\newcommand{\tint}{\text{int}}

\newcommand{\bx}{\mathbf{x}}
\newcommand{\bX}{\mathbf{X}}

\newcommand{\sgn}{\text{sgn}}

\newcommand{\dsep}{\text{DegNSEP}^*}
\newcommand{\pertsep}{\text{PertSEP}^*}
\newcommand{\pertnsep}{\text{PertNSEP}^*}
\newcommand{\dnsep}{\text{DegSEP}^*}

\newcommand{\tg}{\tilde\nabla}

\newcommand{\bbR}{\mathbb{R}}

\newcommand{\bbE}{\mathbb{E}}
\newcommand{\bbP}{\mathbb{P}}

\newcommand{\calD}{\mathcal{D}}
\newcommand{\LRD}{\mathrm{LR}_\calD}
\newcommand{\supp}{\mathrm{supp}}
\newcommand{\dsepD}{\mathrm{DegNSEP}^*_\calD}

\newcommand{\Null}{\text{null}}

\newcommand{\tr}{\mathrm{Tr}}

\DeclareMathOperator*{\argmax}{arg\,max}

\usepackage{fancyhdr} 
\usepackage[nottoc,notlof,notlot]{tocbibind} 
\usepackage[titles,subfigure]{tocloft} 

\title{Condition Number Analysis of Logistic Regression, and its Implications for Standard First-Order Solution Methods }
\author{Robert M. Freund\thanks{MIT Sloan School of Management, 77 Massachusetts Avenue, Cambridge, MA   02139
({mailto:  rfreund@mit.edu}).  This author's research is supported by AFOSR Grant No. FA9550-15-1-0276 and the MIT-Chile-Pontificia Universidad Cat\'olica de Chile Seed Fund, and the MIT-Belgium Universit\'{e} Catholique de Louvain Fund.}
\and Paul Grigas\thanks{Department of Industrial Engineering and Operations Research, University of California, Berkeley, CA 94720
({mailto:  pgrigas@berkeley.edu}).  This author's research is supported by NSF Awards CCF-1755705 and CMMI-1762744.}
\and Rahul Mazumder\thanks{MIT Sloan School of Management, 77 Massachusetts Avenue, Cambridge, MA   02139
({mailto:  rahulmaz@mit.edu}).   This author's research was partially supported by ONR-N000141512342, ONR-N000141812298 (YIP) and NSF-IIS-1718258.}}

\begin{document}
\maketitle

\begin{abstract}
Logistic regression is one of the most popular methods in binary classification, wherein estimation of model parameters is carried out by solving the maximum likelihood (ML) optimization problem, and the ML estimator is defined to be the optimal solution of this problem. It is well known that the ML estimator exists when the data is non-separable, but fails to exist when the data is linearly separable.  First-order methods are the algorithms of choice for solving large-scale instances of the logistic regression problem.  In this paper, we introduce a pair of condition numbers that measure the degree of non-separability or separability of a given dataset in the setting of binary classification, and we study how these condition numbers relate to and inform the properties and the convergence guarantees of first-order methods.  When the training data is non-separable, we show that the degree of non-separability naturally enters the analysis and informs the properties and convergence guarantees of two standard first-order methods:  steepest descent (for any given norm) and stochastic gradient descent.  Expanding on the work of Bach, we also show how the degree of non-separability enters into the analysis of linear convergence of steepest descent (without needing strong convexity), as well as the adaptive convergence of stochastic gradient descent.  When the training data is separable, first-order methods rather curiously have good empirical success -- a behavior that is not well understood in theory.  In the case of separable data, we demonstrate how the degree of separability enters into the analysis of $\ell_2$ steepest descent and stochastic gradient descent for delivering approximate-maximum-margin solutions with associated computational guarantees as well.  This suggests that first-order methods can lead to statistically meaningful solutions in the separable case, even though the ML solution does not exist.
\end{abstract}

\section{Introduction}
Logistic regression is arguably one of the most popular methods for binary classification -- in contrast to SVM-based classifiers, logistic regression provides estimates of the probability of class membership, which is useful for uncertainty quantification and statistical inference.
Moreover, the logistic loss function (and its multiclass extension, the cross-entropy loss) is an essential ingredient of several popular and powerful statistical methods, such as boosting \cite{logitBoost}, kernel methods \cite{zhu2005kernel}, and deep learning \cite{Goodfellow-et-al-2016}.

Let us recall the setting of binary classification and logistic regression. 
Given a binary response $y \in \{-1, 1\}$ and feature vector $\bx \in \bbR^p$, we consider a probability model of the form:
\begin{equation}\label{prob-model}
\bbP\left( y = +1 ~|~ \bx \right) ~=~ \frac{1}{1 + \exp(-\beta^T\bx)} \ ,
\end{equation}
for a vector of coefficients $\beta \in \bbR^p$. The standard procedure for estimating the (unknown) coefficients $\beta$ in \eqref{prob-model} based on a given dataset of $n$ observations $(\bx_1, y_1), \ldots, (\bx_n, y_n)$ is to apply the principle of maximum likelihood (ML) estimation. After some basic algebraic manipulations, maximum likelihood estimation yields the following convex optimization problem: 
\begin{equation}\label{poi-logit}
\begin{array}{rccl}
\mathrm{LR} \ : \ \ \ L_{n}^* := & \min\limits_{\beta} &  & L_{n}(\beta)  \ := \ \frac{1}{n}\sum_{i=1}^n \ln\left(1 + \exp\left(-y_i\beta^T\bx_i\right) \right) \\ \\
& \text{s.t.} & & \beta \in \mathbb{R}^p \ .
\end{array}
\end{equation}
The above objective function $L_n(\cdot)$ is referred to as the logistic loss function, the $i^{\mathrm{th}}$ term of which measures the value of the logistic loss $t \mapsto \ln(1 + \exp(-t))$ on the $i^{\text{th}}$ observation $(\bx_i, y_i)$.

Due in part to classical studies~\cite{silvapulle1981existence,albert1984existence} pertaining to the existence of a ML estimator for the logistic regression problem as well as the prevalence of support vector machines (SVMs), it has become natural and customary to characterize binary classification problems in terms of their \emph{separability} properties.  More formally, a given dataset is either \emph{separable}, in which case the set of observations with $y_i = +1$ may be separated from the set of observations with $y_i = -1$ by a (linear) hyperplane, or is \emph{non-separable}, in which case no such linear separator exists. 
Earlier work in the statistics literature by~\cite{silvapulle1981existence,albert1984existence} have shown that a ML estimator for logistic regression exists when the data is non-separable, and it does not exist when the data is separable. Fairly recently, \cite{candes2018phase} studies phase transitions of the existence of a solution when the features arise from a Gaussian ensemble. A related important theme pertains to algorithms for computing a solution to problem LR. Informally, it is well-known that computational schemes for fitting a logistic regression model by solving the problem LR are ``well-behaved'' when the dataset is non-separable.  
Indeed, by simply examining \eqref{prob-model}, it is evident that if the data is truly generated according to a probabilistic model satisfying \eqref{prob-model}, then with enough samples the dataset will eventually be non-separable; therefore non-separability is somehow an essential characteristic of logistic regression.  On the other hand, separability of a dataset (especially when $n > p$) suggests that there actually is a linear model that can discriminate between the two classes with high accuracy.  In this case, it is well understood that the SVM method may be used to identify a ``good" linear separator.  Furthermore, although the logistic loss function encourages models that linearly separate the data, it does not distinguish between such models, i.e., all linear separators are ``equally favored" by the logistic loss function.  This is in contrast to the SVM method which can be used to find a particularly good (i.e., large margin) linear separator.  The behavior of computational schemes for LR when the dataset is separable is not so well understood in theory, 
though there is recent work on first-order methods \cite{telgarsky2018}, \cite{soudry2017implicit},  \cite{srebro2018a}, \cite{gunasekar2018characterizing}, \cite{srebro2018b}. One of the main goals of this paper is to formalize the ``informal'' computational and statistical intuitions regarding logistic regression and to provide formal results that validate (or run counter to) such intuitive statements.  In particular, a natural set of questions is:  can we quantify the degree of non-separability or separability of a particular dataset, and how might such a formalism inform the computational or statistical properties of solution methods for LR? Herein, we address these questions in the context of first-order methods, which are the methods of choice in the high-dimensional regime ($n \gg 0$ and/or $p \gg 0$).

In recent years, with growing volumes of data, there has been an ever-increasing need to fit accurate logistic regression models to very large datasets with $n \gg 0$ and/or $p \gg 0$. First-order methods for tackling the problem LR are appealing in this large-scale regime for several reasons. First, the computational cost per iteration of first-order methods is relatively low compared to alternatives such as Newton's method (i.e., iteratively reweighted least squares). Second, first-order methods often tend to produce statistically interesting solutions \emph{before} they reach convergence. In particular, several first-order methods such as gradient descent and its generalizations -- steepest descent and stochastic gradient descent -- are known to impart implicit regularization which induces models with good out-of-sample performance on the interior of the sequence of coefficient iterates. Moreover, at least one special case of steepest descent, namely greedy coordinate descent, also imparts desirable sparsity properties along the sequence of coefficient iterates. In this paper, we focus on the method of steepest descent in an arbitrary given norm $\|\cdot\|$ -- a method that encompasses both standard gradient descent and greedy coordinate descent, among others -- as well as the method of stochastic gradient descent (SGD), which is particularly appealing for problems with $n \gg 0$ since SGD only needs to sample a handful of data observations at each iteration.

Towards improving our understanding of steepest descent and SGD for logistic regression, we introduce a pair of condition numbers that measure the degree of non-separability or separability of the dataset $(\bx_1, y_1), \ldots, (\bx_n, y_n)$. In the case when the data is not separable, we introduce a condition number $\dsep$ that precisely quantifies the degree of non-separability of the dataset, namely datasets that are ``more non-separable" have larger values of $\dsep$.  We then show that $\dsep$ naturally informs the computational guarantees of steepest descent in the sense that the guarantees improve when $\dsep$ is larger. Furthermore, we extend the definition of $\dsep$ to measure the degree of non-separability of an arbitrary distribution over the data, which allows us to analyze the role of $\dsep$ in the computational guarantees of SGD in full generality. In particular, we demonstrate that better convergence bounds -- and therefore less data samples from a statistical learning point of view -- are achieved when $\dsep$ is larger.  In the case of separable data, we use the well-known concept of the margin \cite{cover1965geometrical}, which we refer to as $\dnsep$, to precisely quantify the degree of separability of the dataset, whereby datasets that are ``more separable" have larger values of $\dnsep$.  We then develop computational guarantees for both $\ell_2$ steepest descent and SGD that naturally depend on $\dnsep$ and that demonstrate convergence towards approximate-maximum-margin solutions.
We also demonstrate that both $\dsep$ and $\dnsep$ may also be interpreted through the lens of data perturbations and the ``distance to ill-posedness'' introduced by Renegar \cite{Reneg94}.

There has recently been other research activity on the analysis of the performance of steepest descent and SGD for solving the logistic regression optimization problem \eqref{poi-logit}.  Ji and Telgarsky \cite{telgarsky2018}, Soudry et al. \cite{soudry2017implicit}, Nacson et al. \cite{srebro2018a}, Gunasekar et al. \cite{gunasekar2018characterizing}, and Nacson et al. \cite{srebro2018b} analyze the convergence properties of steepest descent and/or SGD in terms of the loss function values and iterate values when the problem instance is separable (or partially separable, this latter case not having received any previous attention that we are aware of).  We discuss the results in these papers relative to ours in the relevant sections herein.

The paper is organized as follows.  In Section \ref{hereiam} we present a pair of condition numbers for logistic regression instances, one for instances that are non-separable and another for instances that are separable.  In Section \ref{imp} we examine the steepest descent algorithm (in any given norm) and show how the degree of non-separability naturally informs computational guarantees of steepest descent.
In the separable case, we develop computational guarantees for $\ell_2$ steepest descent that are informed by the degree of separability and show convergence towards an approximate-maximum-margin solution. 
Expanding on Bach \cite{bach2010, bach2014adaptivity}, we also show how the degree of non-separability enters into the analysis of linear convergence of steepest descent (without needing strong convexity).
In Section \ref{nye} we examine the stochastic gradient descent (SGD) method and we show how our condition numbers inform the computational guarantees of SGD.  In the non-separable case, we show how the degree of non-separability informs standard guarantees of SGD as well as the adaptive guarantee developed by Bach in \cite{bach2014adaptivity}. In the separable case, we develop computational guarantees for SGD that are informed by the degree of separability and show convergence in probability towards an approximate-maximum-margin solution. 

\subsection{Notation}
For a vector $x \in \mathbb{R}^p$, $x_j$ denotes the $j^{\text{th}}$ coordinate; we use superscripts to index vectors in a sequence $\{x^k\}$. Let $e_j$ denote the $j^{\text{th}}$ unit vector in $\mathbb{R}^p$, and let $e = (1, \ldots, 1)$.  We will make use of a generic given norm $\| \cdot \|$ on $\bbR^p$ as well as the $\ell_q$ norm denoted by $\|\cdot\|_q$ with unit ball $B_q$ for $q \in [1, \infty]$. For the given norm $\|\cdot\|$, $\|\cdot\|_\ast$ denotes the dual norm defined by $\|s\|_\ast = \max\limits_{x : \|x\| \leq 1} s^Tx$. Let $\mathrm{Dist}(v, S) := \min_{x \in S} \|x-v\|$ denote the distance from a point $v$ to a set $S$.
Let $\|v\|_0$ denote the number of non-zero coefficients of the vector $v$. 

For $A \in \mathbb{R}^{n \times p}$, let $\|A\|_{\cdot, q} := \max\limits_{x : \|x\| \leq 1} \|Ax\|_{q}$ be the operator norm using the given norm $\|\cdot\|$, and let $\|A\|_{q_1, q_2} := \max\limits_{x : \|x\|_{q_1} \leq 1}\|Ax\|_{q_2}$ be the operator norm where the given norm is the $\ell_{q_1}$ norm. Furthermore, let $\Null(A) := \{x \in \bbR^p : Ax = 0\}$ denote the null space of $A$. For a symmetric matrix $M$, we write ``$M \succeq 0$'' to denote that $M$ is positive semidefinite, ``$M \succ 0$'' to denote that $M$ is positive definite, and let $\lambda_{\min}(M)$ denote the smallest eigenvalue of $M$.

For a scalar $\alpha$, $\sgn(\alpha)$ denotes the sign of $\alpha$, and $\alpha^+$, $\alpha^-$ denote the positive and negative parts of $\alpha$, respectively.  The notation ``$\tilde v \leftarrow \argmax\limits_{v  \in S} \{f(v)\}$'' denotes assigning $\tilde v$ to be any optimal solution of the problem $\max\limits_{v  \in S} \{f(v)\}$.

\section{Logistic Regression, and a Pair of Condition Numbers for Non-Separable and Separable Training Data}\label{hereiam}

Let us review the setting and notation of logistic regression and the basic properties of the optimization problem LR.  
Recall that we have $n$ observed training data points $(\bx_1, y_1)\ldots, (\bx_n, y_n)$ where $\bx_i \in \mathbb{R}^p$ is the vector of feature values and $y_i \in \{-1,1\}$ is the (binary) class of observation $i$, for $i=1, \ldots, n$. Let $\bX$ be the matrix whose $i^{\mathrm{th}}$ row is $\bx_i$, for $i=1, \ldots, n$.  The well-known logistic loss function is $L_{n}(\cdot) : \mathbb{R}^p \to \mathbb{R}$ defined by:
\begin{equation}\label{lloss}
L_{n}(\beta) := \frac{1}{n}\sum_{i=1}^n \ln\left(1 + \exp\left(-y_i\beta^T\bx_i\right) \right) \ ,
\end{equation}
where $\beta \in \bbR^p$.  Throughout the paper, we denote the univariate logistic loss function by $\ell(t):= \ln(1+ \exp(-t))$, the gradient of $L_n(\cdot)$ by $\nabla L_n(\cdot)$, and the Hessian of $L_n(\beta)$ by $H(\beta)$. 

As mentioned previously, LR, the problem of minimizing the logistic loss function $L_{n}(\cdot)$ over $\beta \in \bbR^p$ arises from maximum likelihood estimation for the model \eqref{prob-model}.  
Note that $L_{n}(\beta) > 0$ for all $\beta$, hence $L_{n}^* \geq 0$ and is therefore finite.  Unlike, for example, the least-squares loss function for linear regression, it is not clear {\it a priori} if the logistic regression problem LR has an optimal solution.  Indeed, in the case when the data is separable, i.e., there exists a vector $\beta \in \mathbb{R}^p$ satisfying $y_i \beta^T\bx_i> 0$ for $i=1,\ldots,n$, then $L_{n}(\theta \beta) \rightarrow 0=L_{n}^* $ as $\theta \rightarrow +\infty$, and hence LR does not attain its optimum.

In order to better understand the behavior of the logistic regression problem in general as well as the behavior of the steepest descent and SGD for logistic regression, we now develop a pair of condition numbers that measure the degree of non-separability and separability of the dataset $(\bx_1, y_1)\ldots, (\bx_n, y_n)$.  In particular, we show in this section that the behavior of LR in terms of existence of optima, as well as the existence of low-norm optima, can be characterized in terms of these condition numbers.

\subsection{Non-Separable Data}\label{nonseparable}

Let us first consider the case of non-separable data.  We say that observation $i$ is {\em correctly classified} by $\beta$ if $y_i \beta^T\bx_i > 0$, and is {\em misclassified} by $\beta$ if $y_i \beta^T\bx_i \le 0$.  Letting $Y$ denote the diagonal matrix whose $i^{\mathrm{th}}$ component is $y_i$, then with this notation observation $i$ is correctly classified or misclassified by $\beta$ if $(Y\bX\beta)_i >0$ or $(Y\bX\beta)_i \le 0$, respectively.  We say that the training data is {\em non-separable} if there is no $\beta$ that correctly classifies every observation, i.e., there is no $\beta$ that satisfies $Y\bX\beta >0$, and in this case we write ``$(\bX, y)$ is not separable'' to denote that the data $(\bX,y)$ are not separable.

Clearly, some non-separable datasets might be ``more non-separable'' than others, so let us now introduce a way to measure the extent to which the dataset is non-separable.  Let $\|\cdot\|$ denote the given norm on the space $\mathbb{R}^p$ of model coefficients $\beta$.  We define the degree of non-separability of the training data (with respect to the norm $\|\cdot\|$) to be:

\begin{equation}\label{dsep}
\begin{array}{rccl}
 \ \ \ \mathrm{DegNSEP}^* := & \displaystyle\min_{\beta \in \mathbb{R}^p} &  & \frac{1}{n}\sum_{i=1}^n [y_i\beta^T\bx_i]^{-} \\ \\
& \text{s.t.} & & \|\beta\| = 1 \ ,
\end{array}
\end{equation}

which states that $\dsep$ is the smallest (over all normalized models $\beta$) average misclassification error of the model $\beta$ over the $n$ observations.  We emphasize that here and for the remainder of this paper that the norm $\|\cdot\|$ on the space of model coefficients $\mathbb{R}^p$ is given and fixed. (In Sections \ref{sdsubsection_sep} and \ref{nye} we will require the norm $\|\cdot\|$ to be the $\ell_2$ norm, but otherwise it is generic.) 

Define $\beta^0 := 0 \in \mathbb{R}^p$.  Noticing that $L_n(\beta^0) = \ln(2)$, the following object measures the maximum distance from any $\beta$ in the level set of $\beta^0$ -- namely $\{ \beta \in \mathbb{R}^p : L_n(\beta) \le \ln(2)\}$ -- to the set of optimal solutions of LR:

\begin{equation}\label{long} \mathrm{Dist}_0 = \max\limits_{\beta : L_{n}(\beta) \le \ln(2)}  \  \left\{ \min\limits_{\beta^* : L_{n}(\beta^*) = L_{n}^*}\|\beta - \beta^*\| \right\} \ . \end{equation}

The following proposition shows that the behavior of the logistic regression problem LR can be characterized in terms of the degree of non-separability $\dsep$.\medskip

\begin{proposition}\label{strictly2}  If $\mathrm{DegNSEP}^* > 0$, then:
\begin{enumerate}
\item[(i)] there is a unique optimal solution $\beta^\ast$ of the logistic regression problem LR,
\item[(ii)] $H(\beta^\ast) \succ 0$ ,
\item[(iii)] $\|\beta^*\| ~\le~ \displaystyle\frac{L_n^*}{\mathrm{DegNSEP}^*} ~\le~  \displaystyle\frac{\ln(2)}{\mathrm{DegNSEP}^*}$ , and
\item[(iv)] $\mathrm{Dist}_0 ~\le~ \displaystyle\frac{\ln(2)+L_n^*}{\mathrm{DegNSEP}^*} ~\le~ \displaystyle\frac{2\ln(2)}{\mathrm{DegNSEP}^*}$ .
\end{enumerate}
\end{proposition}
\noindent {\bf Proof:} Suppose that $\dsep > 0$. Notice that this implies that $\Null(\bX) = \{0\}$ since if this is not the case, then there exists $\bar\beta \in \Null(\bX)$ with $\|\bar\beta\| = 1$ which implies that $\dsep = 0$.
For any $\beta \in \bbR^p$, a simple calculation yields that $H(\beta) = \tfrac{1}{n}\bX^TG\bX$ where $G$ is the $n \times n$ diagonal matrix whose $i^{\mathrm{th}}$ component is $G_{ii} = \ell^{\prime\prime}(y_i\beta^T\bx_i) = \frac{\exp(y_i\beta^T\bx_i)}{(\exp(y_i\beta^T\bx_i) + 1)^2} > 0$. Therefore, for any $\beta \in \bbR^p$, we have that $\Null(H(\beta)) = \Null(\bX) = \{0\}$, and it then follows that $H(\beta) \succ 0$.  This implies that $L_n(\cdot)$ is globally strictly convex.   

Notice that $\ln(1+e^{-t}) \ge t^-$ for any $t$, and hence the objective function of \eqref{dsep} satisfies $L_n(\beta) \ge \frac{1}{n}\sum_{i=1}^n [y_i\beta^T\bx_i]^{-} $ for any $\beta$.  It then follows from \eqref{dsep} that $L_n(\beta) \ge \dsep \|\beta\| $ for all $\beta \in \mathbb{R}^p$, which rearranges to:

\begin{equation}\label{january4} \|\beta\| \le \frac{L_n(\beta)}{\dsep} \ \ \ \ \ \mbox{for~all~} \beta \in\mathbb{R}^p \ . \end{equation}

Notice that $L_n(0)=\ln(2)$, and therefore the level set $\{ \beta \in \mathbb{R}^p : L_n(\beta) \le \ln(2)\} \subset \{ \beta \in \mathbb{R}^p : \|\beta\| \le \ln(2)/\dsep \}$ and hence is a nonempty compact set.  
It then follows from the continuity of $L_n(\cdot)$ in conjunction with the Weierstrass Theorem that LR attains its optimum. Since $L_n(\cdot)$ is strictly convex there is a unique optimal solution $\beta^\ast$, which proves {\em (i)}. By the previous discussion we have that $H(\beta^\ast) \succ 0$, which proves {\em (ii)}, and it follows from \eqref{january4} that $\|\beta^*\| \le \frac{L_n^*}{\dsep}$, which proves {\em (iii)}.  If $\beta$ satisfies $L_n(\beta) \le \ln(2)$ it follows that
$$\|\beta - \beta^*\| \le \|\beta\| + \|\beta^*\| \le \frac{\ln(2)}{\dsep}  + \frac{L_n^*}{\dsep}  \ , $$
which then implies {\em (iv)}. \qed

Part {\em (iii)} of Proposition \ref{strictly2} states that the norm of the unique optimal solution of the logistic regression problem LR is bounded inversely proportional to $\dsep$, and part {\em (iv)} of Proposition \ref{strictly2} presents a similar bound on  $\mathrm{Dist}_0$. Part {\em (ii)} of Proposition \ref{strictly2} states that $H(\beta^\ast)$ is positive definite, i.e., that the logistic loss function is locally strongly convex at the optimum $\beta^\ast$. We can measure the degree of local strong convexity by defining, for any symmetric positive semidefinite matrix $M$, the local strong convexity constant of $M$ (with respect to the norm $\|\cdot\|$) to be:
\begin{equation}\label{local_strong_constant}
\begin{array}{rccl}
 \ \ \ \nu^\ast(M) := & \displaystyle\min_{\beta \in \mathbb{R}^p} &  & \beta^TM\beta \\ \\
& \text{s.t.} & & \|\beta\| = 1 \ .
\end{array}
\end{equation}
Part {\em (ii)} of Proposition \ref{strictly2} immediately implies that $\nu^\ast(H(\beta^\ast)) > 0$. Notice that when the norm $\|\cdot\|$ is the $\ell_2$ norm, then $\nu^\ast(M) = \lambda_{\min}(M)$.
It also follows from norm equivalence that there exist constants $C_1$ and $C_2$ with $0 < C_1 \leq C_2$ such that $C_1\lambda_{\min}(M) \leq \nu^\ast(M) \leq C_2\lambda_{\min}(M)$ for all symmetric positive semidefinite matrices $M$.
Proposition \ref{local_strong_prop} below provides a lower bound on $\nu^\ast(H(\beta^\ast))$ that depends entirely on $\dsep$ and magnitude properties of $\bX$. 

\begin{proposition}\label{local_strong_prop}
If $\mathrm{DegNSEP}^* > 0$, then:
$$\nu^\ast(H(\beta^\ast)) ~\geq~ \tfrac{1}{4n}\nu^\ast(\bX^T\bX)\exp \left(-\frac{\ln(2)\|\bX\|_{\cdot,\infty}}{\mathrm{DegNSEP}^*}\right) > 0 \ ,$$ where $\beta^\ast$ is the unique optimal solution of LR.
\end{proposition}
\noindent {\bf Proof:}
Recall that $H(\beta^\ast) = \tfrac{1}{n}\bX^TG\bX$ where $G$ is the $n \times n$ diagonal matrix whose $i^{\mathrm{th}}$ component $G_{ii} = \ell^{\prime\prime}(y_i(\beta^\ast)^T\bx_i)$ satisfies $$ G_{ii} = \frac{\exp(y_i(\beta^\ast)^T\bx_i)}{(\exp(y_i(\beta^\ast)^T\bx_i) + 1)^2} \geq \tfrac{1}{4}\exp(-|y_i(\beta^\ast)^T\bx_i|) \geq \tfrac{1}{4}\exp(-\|\bx_i\|_\ast\|\beta^\ast\|) \geq \tfrac{1}{4}\exp\left(-\tfrac{\ln(2)\|\bX\|_{\cdot,\infty}}{\mathrm{DegNSEP}^*}\right) \ , $$
where the final inequality uses part {\em (iii)} of Proposition \ref{strictly2} and $\|\bX\|_{\cdot, \infty} = \max\limits_{i \in \{1, \ldots, n\}} \|\bx_i\|_\ast$. Therefore, for any $\beta \in \bbR^p$, we have
\begin{equation*}
\beta^T H(\beta^\ast)\beta = \tfrac{1}{n}(\bX\beta)^TG(\bX\beta) \geq \tfrac{1}{4n}\exp\left(-\tfrac{\ln(2)\|\bX\|_{\cdot,\infty}}{\mathrm{DegNSEP}^*}\right)\beta^T(\bX^T\bX)\beta \ ,
\end{equation*}
and taking the minimum of both sides of the above inequality over all $\beta$ satisfying $\|\beta\| = 1$ yields the desired result.
\qed

It turns out that $\dsep$ can also be interpreted as the minimal {\em perturbation} of the data for which the perturbed problem data is separable.  This will be given a precise definition and proof in Section \ref{secperturbation}.

\subsection{Separable Data}\label{separable}

We say that the training data is {\em separable} if there exists $\beta$ for which $Y\bX \beta > 0$, i.e., there is a model $\beta$ that correctly classifies every observation, and we write ``$(\bX, y)$ is separable'' to denote that the data $(\bX,y)$ are separable.  Akin to the case of non-separable data, some separable datasets might be ``more separable'' than others.  Employing the standard lens of statistical machine learning \cite{ESLBook}, we can measure the degree of separability using the well-known concept of the ``margin'' \cite{cover1965geometrical}, which we now review for completeness.  Suppose that $(\bX, y)$ is separable, and let $\beta$ be a model that correctly classifies all observations, namely $Y\bX \beta >0$.  Then the {\em margin} of $\beta$ is denoted by $\rho(\beta)$ and is defined to be the least classification value of $\beta$ over all observations: $$\rho(\beta) :=  \displaystyle\min_{i \in \{1, \ldots, n \}} [ y_i\beta^T\bx_i ] \ . $$  We define the degree of separability of the data to be the maximum margin over all normalized models $\beta$, namely:

\begin{equation}\label{dnsep}
\begin{array}{rll}
 \ \ \ \dnsep := & \displaystyle\max_{\beta}   & \rho(\beta)\\ \\
& \text{s.t.}  & \|\beta\| \le 1 \ . \end{array}
\end{equation}\medskip

\begin{proposition}\label{strictly3}  If $(\bX, y)$ is separable, then $\mathrm{DegSEP}^* >0$, $L_n^* = 0$, and LR does not attain its optimum.
\end{proposition}

\noindent {\bf Proof:} If $(\bX, y)$ is separable, it follows from the definition of the margin that $\mathrm{DegSEP}^* >0$, and there exists a vector $\bar\beta \in \mathbb{R}^p$ satisfying $\|\bar\beta\| \le  1$ and $y_i \bar\beta^T\bx_i \ge \dnsep > 0$ for $i=1,\ldots,n$, whereby $L_{n}(\theta \bar\beta) \rightarrow 0$ as $\theta \rightarrow +\infty$.  Since $L_n(\beta) >0$ for any $\beta$, it also follows that $L_n^* \ge 0$, and hence $L_n^* =0$ and LR does not attain its optimum.\qed

The following lemma relates the margin function $\rho(\cdot)$ to the gradient of the logistic loss function.  In Lemma \ref{margin-lemma} and elsewhere in this paper we use the convention $\ln(a) = -\infty$ for $a \le 0$, i.e., $\ln(\cdot)$ is an extended-real-valued concave function.

\begin{lemma}\label{margin-lemma}
Suppose that the data $(\bX, y)$ is separable, i.e., $\mathrm{DegSEP}^* > 0$. Then for any $\beta \in \bbR^p$ it holds that:
\begin{equation*}
\rho(\beta) ~\geq~ \ln\left(\tfrac{\mathrm{DegSEP}^*}{n\|\nabla L_n(\beta)\|_\ast} - 1\right) \ .
\end{equation*}
\end{lemma}
\noindent {\bf Proof:}
Recall that 
\begin{equation*}
\rho(\beta) :=  \displaystyle\min_{i \in \{1, \ldots, n \}}[y_i\beta^T\bx_i] = \displaystyle\min_{i \in \{1, \ldots, n \}} (Y\bX\beta)_i = \displaystyle\min_{w \in \Delta_n} w^TY\bX\beta \ ,
\end{equation*} where $\Delta_n := \{ w \in \mathbb{R}^n: e^Tw=1, \ w \ge 0\}$.  By minimax strong duality it holds that:
\begin{equation}\label{dual-rho}
\dnsep :=  \displaystyle\max_{\beta: \|\beta\| \le 1}   \  \rho(\beta) \ \ = \ \   \displaystyle\min_{w \in \Delta_n} \|\bX^T Y w\|_* \ .
\end{equation}
Straightforward calculation yields that the gradient of the logistic loss function satisfies $\nabla L_n(\beta) = \tfrac{1}{n}\bX^TYw^\ast(\beta)$ where $w^\ast(\beta) \in \bbR^n$ is the vector defined component-wise by:
\begin{equation}\label{w-def}
w^\ast(\beta)_i = \frac{1}{1+\exp(y_i \beta^T\bx_i)} \ \ , \ \  i =1, \ldots, n \ .
\end{equation}
Since $w^\ast(\beta) > 0$ it follows that $w^\ast(\beta)/\|w^\ast(\beta)\|_1$ is a feasible solution of the right-most optimization problem in \eqref{dual-rho}, which implies that $\frac{\|\bX^TYw^\ast(\beta)\|_\ast}{\|w^\ast(\beta)\|_1} \geq \dnsep$.
Thus, since $\dnsep > 0$, it holds that
\begin{equation*}
\|w^\ast(\beta)\|_1 \leq \frac{\|\bX^TYw^\ast(\beta)\|_\ast}{\dnsep} = \frac{n\|\nabla L_n(\beta)\|_\ast}{\dnsep} \ .
\end{equation*}
In particular, for each $i =1, \ldots, n$, it holds that $w^\ast(\beta)_i \leq \tfrac{n\|\nabla L_n(\beta)\|_\ast}{\dnsep}$, which, after simple arithmetic manipulation using \eqref{w-def}, is equivalent to:
\begin{equation}\label{individual-margin-bound}
y_i \beta^T\bx_i ~\geq~ \ln\left(\tfrac{\dnsep}{n\|\nabla L_n(\beta)\|_\ast} - 1\right) \ ,  \  i =1, \ldots, n \ .
\end{equation}
Since \eqref{individual-margin-bound} holds for all $i =1, \ldots, n$, the result is proved.\qed

It turns out that $\dnsep$ can also be interpreted as the minimal perturbation of the data for which the perturbed problem data is non-separable.  This is developed in the following subsection.

\subsection{Data-Perturbation Interpretations of $\dsep$ and $\dnsep$}\label{secperturbation}

In this subsection we show that both $\dsep$ and $\dnsep$ can be interpreted through the lens of data perturbations that alter the status of the dataset -- from non-separable to separable, or {\em vice versa}.  Let us view the feature data matrix $\bX$ as a linear operator, and recall the operator norm notation $\|\bX\|_{\cdot,q} := \max\limits_{\beta : \|\beta\| \leq 1}\|\bX \beta \|_{q}$ on the space $\mathbb{R}^{n \times p}$.

Define:

\begin{equation}\label{evelyn}
\begin{array}{rll}
 \ \ \ \pertsep := & \displaystyle\inf_{\Delta \bX}   & \tfrac{1}{n}\|\Delta \bX\|_{\cdot,1} \\ \\
& \text{s.t.}  & (\bX + \Delta \bX, y) ~\mbox{is~separable} \ . \end{array}
\end{equation}

Then $\pertsep$ is the smallest (or more precisely, the infimum thereof) perturbation $\Delta \bX$ of the feature data $\bX$ which will render the perturbed problem instance $(\bX + \Delta \bX, y)$ separable.  Here the size of the perturbation is measured using the (scaled) operator norm $\tfrac{1}{n}\| \cdot \|_{\cdot,1}$.  Clearly $\pertsep = 0$ if $(\bX, y)$ is separable.  If $\pertsep >0$, then $(\bX, y)$ is not separable; and the smaller $\pertsep$ is, the closer the data is to being separable.  Notice that $\|\Delta \bX\|_{\cdot,1}$ scales proportional to $n$, which is counteracted by dividing by $n$ in the objective function of \eqref{evelyn}.  The following result shows that the condition number $\dsep$ introduced and used in Section \ref{nonseparable} can be alternatively interpreted as the smallest data perturbation for which the perturbed data  $(\bX+\Delta \bX, y)$ is separable.\medskip

\begin{proposition}\label{jimrenegar1} For any dataset $(\bX, y)$ it holds that $\dsep = \pertsep$.  
\end{proposition}

{\bf Proof:}  See Appendix \ref{June29}. \qed \medskip

Let us also define:

\begin{equation}\label{drew}
\begin{array}{rll}
 \ \ \ \pertnsep := & \displaystyle\inf_{\Delta \bX}   & \|\Delta \bX\|_{\cdot,\infty} \\ \\
& \text{s.t.}  & (\bX + \Delta \bX, y) ~\mbox{is~non-separable} \ . \end{array}
\end{equation}

Then $\pertnsep$ is the smallest perturbation $\Delta \bX$ of the feature data $\bX$ which will render the perturbed problem instance $(\bX + \Delta \bX, y)$ non-separable.  Here the size of the perturbation is measured using the operator norm $\| \cdot \|_{\cdot,\infty}$.  Clearly $\pertnsep = 0 $ if $(\bX, y)$ is not separable.  If $\pertnsep >0$, then $(\bX, y)$ is separable; and the smaller $\pertnsep$ is, the closer the data is to being non-separable.  Notice that we use a different operator norm in the definition of $\pertnsep$ than that used in the definition of $\pertsep$.  The following result shows that the condition number $\dnsep$ introduced and used in Section \ref{separable} can be alternatively interpreted as the smallest data perturbation for which the perturbed data  $(\bX+\Delta \bX, y)$ is non-separable.\medskip

\begin{proposition}\label{jimrenegar2} For any dataset $(\bX, y)$ it holds that $\dnsep = \pertnsep$.  
\end{proposition}

{\bf Proof:}  See Appendix \ref{June29} as well. \qed \medskip

Any given dataset $(\bX, y)$ is either non-separable (in which case $\dnsep = 0$) or is separable (in which case $\dsep = 0$).   Borrowing from the lexicon of Renegar \cite{Reneg94}, we may call the dataset $(\bX, y)$ ``ill-posed'' if both $\dsep = 0 $ and $\dnsep = 0$, in which case the dataset is non-separable but there is an arbitrarily small perturbation of the data that renders the perturbed dataset separable.  It can easily be verified that an example of such a dataset is:
$$ \bX = \begin{pmatrix*}[r] 1 & 0 & -1 \cr 0 & -1 & 1 \cr -1 & -2 & 3 \cr 2 & 1 & -3 \end{pmatrix*} \ , \ \ y = \begin{pmatrix*}[r] 1  \cr -1  \cr 1 \cr -1 \end{pmatrix*} \ . $$

\section{Informing Standard Deterministic First-Order Solution Methods for Logistic Regression}\label{imp}

In this section we show how the two condition numbers $\dsep$ and $\dnsep$ inform the computational properties and guarantees of a standard deterministic first-order solution method for logistic regression, namely the steepest descent method in any given norm $\|\cdot\|$.  After briefly reviewing steepest descent in the setting of smooth convex optimization in Section \ref{sdwithf}, we examine steepest descent as applied to logistic regression in two cases: in Section \ref{sdsubsection} we consider the case of non-separable data and examine steepest descent for any given norm, and in Section \ref{sdsubsection_sep} we consider the case of separable data and examine steepest descent for the $\ell_2$ norm.

\subsection{Brief Review of Steepest Descent}\label{sdwithf}
We review the steepest descent method for solving the unconstrained optimization problem:
\begin{equation}\label{poi4}
\mbox{(P):} \ \ f^*:=  \ \min\limits_{x \in \mathbb{R}^p} \ f(x) \ ,
\end{equation}
where $f(\cdot) : \mathbb{R}^p \to \mathbb{R}$ is a differentiable convex function, and we assume that $f^*$ is finite but it is not necessarily attained.  Let $\|\cdot\|$ be the norm on the variables $x \in \mathbb{R}^p$.  At a given iterate $\bar x$, the steepest descent direction is defined to be the negative of the normalized direction $\bar d$ that maximizes $\nabla f(\bar x)^T d$, namely:  $\bar d \gets \arg\max_d \{\nabla f(\bar x)^T d : \|d \| \le 1 \}$.  The formal statement of steepest descent in the norm $\|\cdot\|$ is presented in Algorithm \ref{sdm}.\medskip 

\floatname{algorithm}{Algorithm}
\begin{algorithm}
\caption{Steepest Descent in the norm $\|\cdot\|$}\label{sdm}
\begin{algorithmic}
\STATE Initialize at $x^0 \in \mathbb{R}^p$, $k \leftarrow 0$\medskip

At iteration $k$:
\STATE 1. Compute $\nabla f(x^k)$
\STATE 2. Compute $d^k \gets \arg\max_d \{\nabla f(x^k)^T d : \|d \| \le 1 \}$
\STATE 3. Choose $\alpha_k \geq 0$ and set:
\begin{description}
\item $\ \  x^{k+1} \gets x^k - \alpha_k \cdot d^k$
\end{description}
\end{algorithmic}
\end{algorithm}

To the best of our knowledge, there is no general computational guarantee associated with steepest descent for an arbitrary norm without additional assumptions such as bounded optimal solutions, strong convexity of the function, and/or maximum distances of the starting points and/or the iterates from the set of optimal solutions (or the set of near-optimal solutions).  In the Euclidean norm case ($\|\cdot\| = \|\cdot\|_2$), steepest descent specifies to the classical gradient descent algorithm, see \cite{polyak} and \cite{nesterovBook}.  In the particular case when $\|\cdot\|$ is the $\ell_1$ norm $\|\cdot\|_1$, steepest descent specifies to the greedy coordinate descent method.  For works related to greedy coordinate descent, including block-coordinate descent, cyclic and randomized coordinate descent and variations thereof, see for example Nesterov \cite{Nesterovgcd}, Richt\'{a}rik and Tak\'{a}\v{c} \cite{richgcd}, Schmidt and Friedlander \cite{schmidt2012}, and Beck and Tetruashvili \cite{beckgcd}.

Recall that $f(\cdot)$ is $L$-smooth with respect to the norm $\|\cdot\|$ if $f(\cdot)$ is differentiable and satisfies:
\begin{equation}\label{smoothness} \|\nabla f(x) - \nabla f(y)\|_* \le L \|x-y\| \ \ \ \mbox{for~all~} x,y \in \mathbb{R}^p \ , \end{equation} where $\|\cdot\|_\ast$ is the dual norm of $\|\cdot\|$.

Let ${\cal S}_0$ denote the level set of the initial iterate $x^0$, namely ${\cal S}_0:= \{x \in \mathbb{R}^p : f(x) \le f(x^0)\}$, and let ${\cal S}^*$ denote the set of optimal solutions of \eqref{poi4}, i.e., ${\cal S}^*:= \{x \in \mathbb{R}^p : f(x) =f^*\}$.  Then let $\mathrm{Dist}_0$ denote the largest distance of points in ${\cal S}_0$ to the set of optimal solutions ${\cal S}^*$:
\begin{equation}\label{saturday} \mathrm{Dist}_0 := \max\limits_{x \in {\cal S}_0} \min\limits_{x^* \in {\cal S}^*} \|x-x^*\| \ .
\end{equation}

The following computational guarantees for steepest descent are an amended and extended version of mostly well-known results about traditional gradient descent and greedy coordinate descent, see for example Nesterov \cite{nesterovBook} and Beck and Tetruashvili \cite{beckgcd}.  \medskip

\begin{theorem}\label{sdmcomplexity} {\bf (Computational Guarantees for Steepest Descent in the norm $\|\cdot\|$)} Let $\{x^k\}$ be generated according to the steepest descent method (Algorithm \ref{sdm}) using the step-size sequence $\{\alpha_k \}$ chosen using the ``greedy'' rule:
\begin{equation}\label{sdmquadrule} \alpha_k = \frac{\|\nabla f(x^k)\|_*}{L} \  \  \ \mbox{for~all~} k \ge 0 \ . \end{equation} If $x^0 \notin \cal S^*$, then it holds for all $k \ge 0$ that:

\begin{itemize}
\item[(i)] (optimality gap):
$
f(x^k) - f^*  \ \ \leq \ \ \displaystyle\frac{2L(\mathrm{Dist}_0)^2}{\hat K^0 + k} \ \ < \ \ \frac{2L(\mathrm{Dist}_0)^2}{k}$

\item[] and
\item[(ii)] (gradient bound I):
$ \|\nabla f(x^k)\|_*  \ \ \leq \ \ \sqrt{2L(f(x^k)-f^*)}  \ \ \leq \ \ \displaystyle\frac{2L\mathrm{Dist}_0}{\sqrt{\hat K^0 + k}}$
\end{itemize}
where $\hat K^0 := \displaystyle\frac{2L(\mathrm{Dist}_0)^2}{f(x^0)-f^*}$ .  Furthermore, 
\begin{itemize}
\item[(iii)] (norm bound):
$\|x^k - x^0\| \ \le \ \sqrt{k}\sqrt{\displaystyle\frac{2(f(x^0)-f^*)}{L}}$ , and
\item[(iv)] (gradient bound II):
there exists $i \le k$ for which
$\|\nabla f(x^i)\|_*  \ \ \leq \ \ \sqrt{\displaystyle\frac{2L(f(x^0)-f^*)}{k+1}}$ ,
\end{itemize} where the two inequalities in (i) and the second inequality in (ii) are only relevant when $\mathrm{Dist}_0$ is finite. \qed
\end{theorem}\medskip

For completeness, a self-contained proof of Theorem \ref{sdmcomplexity} is given in Appendix \ref{january2}.

\begin{remark}\label{concordlib} If an exact line-search is used instead of the step-size rule \eqref{sdmquadrule}, then all of the results in Theorem \ref{sdmcomplexity} remain valid except for the norm bound in item {\em (iii)}.  This follows easily from the structure of the proof in Appendix \ref{january2}.\end{remark}

\begin{remark}\label{concordlib2}  In the case of the $\ell_2$ norm, $\mathrm{Dist}_0$ can be replaced by the typically much small quantity $\mathrm{Dist}(x^0, {\cal S}^*)$ in items {\em (i)} and {\em (ii)} of the Theorem \ref{sdmcomplexity}.\end{remark}

\subsection{Informing Steepest Descent for Solving Logistic Regression in the Non-Separable Case}\label{sdsubsection}
Here we show how the condition numbers $\dsep$ and $\dnsep$ inform computational guarantees for steepest descent applied to the logistic regression optimization problem \eqref{poi-logit}.  In order to apply Theorem \ref{sdmcomplexity} to the setting of logistic regression, we use the following smoothness property of the logistic loss function.

\begin{proposition}{\bf{(Lipschitz smoothness of the logistic loss function)}}\label{logistic-props}
The logistic loss function $L_{n}(\cdot)$ is $L=\tfrac{1}{4n}\|\bX\|^2_{\cdot,2}$-smooth with respect to the given norm $\|\cdot\|$ on $\mathbb{R}^p$. \end{proposition}

{\bf Proof:}  See Appendix \ref{tired}. \qed\medskip

\begin{theorem}\label{LogitBoost-complexity1}{\bf (Steepest Descent for Logistic Regression in the norm $\|\cdot\|$:  Non-Separable Case)}
Suppose that steepest descent (Algorithm \ref{sdm}) is initialized at $\beta^0 := 0$, and is implemented using the step-size rule:
\begin{equation}\label{stairmaster1} \alpha_k := \frac{4n\|\nabla L_n(\beta^k)\|_* }{\|\bX\|^2_{\cdot,2}}\ \ \ \mbox{for~all~} k \ge 0 \ .
\end{equation}
When the data is non-separable, for all $k \ge 0$ it holds that:
\begin{itemize}

\item[(i)] (training error):
$$L_{n}(\beta^k)-L^*_n \ \le \ \frac{1}{\displaystyle\frac{1}{\ln(2) - L^*_n} + \displaystyle\frac{k \cdot n \cdot (\mathrm{DegNSEP}^*)^2 }{2\|\bX\|_{\cdot,2}^2(\ln(2))^2} }  \ < \  \displaystyle\frac{2\|\bX\|_{\cdot,2}^2(\ln(2))^2}{k \cdot n \cdot (\mathrm{DegNSEP}^*)^2  } \ , $$

\item[(ii)] (shrinkage): $\|\beta^k\| ~\leq~ \sqrt{k}\left(\frac{1}{\|\bX\|_{\cdot,2}}\right)\sqrt{8n(\ln(2) - L_n^\ast)} ~\leq~ \sqrt{k}\left(\frac{1}{\|\bX\|_{\cdot,2}}\right)\sqrt{8n\ln(2)}$ , and

\item[(iii)] (gradient bound): $\|\nabla L_n(\beta^k)\|_*  \ \leq \ \|\bX\|_{\cdot,2}\sqrt{\frac{(L_n(\beta^k) - L^*_n)}{2n}} \ \leq \ \displaystyle\frac{\|\bX\|_{\cdot,2}^2\ln(2)}{\sqrt{k} \cdot n \cdot \mathrm{DegNSEP}^*}$.

\end{itemize}
\end{theorem}\medskip

\noindent {\bf Proof:} These results are a straightforward application of Theorem \ref{sdmcomplexity}, Proposition \ref{logistic-props}, and Proposition \ref{strictly2}.  Parts {\em (i)} and {\em (ii)} of Theorem \ref{sdmcomplexity} present computational guarantees for the steepest descent method for the step-size rule \eqref{sdmquadrule} in terms of the initial objective function value (which in this case is $L_n(\beta^0) = L_n(0) = \ln(2)$), the Lipschitz constant $L$, and the distance measure $\mathrm{Dist}_0$ defined in \eqref{saturday}. From Proposition \ref{logistic-props} we can take the Lipschitz constant $L$ of the gradient of the logistic loss function $L_n(\cdot)$ to be $L = \tfrac{1}{4n}\|\bX\|^2_{\cdot,2}$.  And from Proposition \ref{strictly2} we have that $\mathrm{Dist}_0 \le \frac{2\ln(2)}{\mathrm{DegNSEP}^*}$ in the case when the data is non-separable.  Substituting these values into the step-size formula \eqref{sdmquadrule} and utilizing parts {\em (i)} and {\em (ii)} of Theorem \ref{sdmcomplexity} yields precisely the step-size rule \eqref{stairmaster1} and the computational guarantees in parts {\em (i)} and {\em (iii)} of the present theorem.  Also, substituting the bound on $L$ into part {\em (iii)} of Theorem \ref{sdmcomplexity} yields part {\em (ii)} of the present theorem. \qed

Notice the manner in which $\dsep$ informs the computational guarantees in Theorem \ref{LogitBoost-complexity1}.  The training error bound scales like $O(1/(\dsep)^2)$, and so the computational guarantee on the training error is smaller for datasets with larger values of $\dsep$.  Also note that the training error bound is invariant under constant re-scaling of the data -- since rescaling all observations by a constant $\gamma$ will rescale both $\dsep$ and $\|\bX\|_{\cdot,2}$ by $\gamma$ and so their quotient is unaffected.  Similarly, the computational guarantee on the norm of the gradient scales like $O(1/\dsep)$, and is smaller for datasets with larger values of $\dsep$ as well.

\subsection{Linear Convergence Results}

In a series of papers, Bach \cite{bach2010, bach2014adaptivity}, as well as Bach and Moulines \cite{bach2013non}, identified a generalized self-concordance property of the logistic loss function that has proven to be useful in analyzing the statistical and computational properties of empirical risk minimization as well as stochastic gradient descent in this setting. In particular, in the case of non-separable data, Bach demonstrates in \cite{bach2014adaptivity} that (averaged) stochastic gradient descent is adaptive to the unknown local strong convexity of the logistic loss function at the optimum. That is, the convergence rate can be improved from $O(1/\sqrt{k})$ to $O(1/\mu k)$ where $\mu$ is the smallest eigenvalue of the Hessian at the optimum. Since it is known that steepest descent achieves linear convergence in the case of a strongly convex objective function, it is natural to ask whether steepest descent is also adaptive to the unknown local strong convexity of the logistic loss function?  Here we answer this question in the affirmative.

We show in the case of non-separable data that steepest descent achieves linear convergence with a rate of linear convergence that naturally depends on the condition number $\dsep$ as well as a measure of local strong convexity at the optimum. Our results provide linear convergence guarantees both in terms of the training error gap as well as the distance to the optimal solution $\beta^\ast$ measured in the given norm.
Moreover, we show that after a certain number of iterations that scales like $O(1/(\mathrm{DegNSEP}^*)^2)$, the rate of linear convergence improves to a faster rate that is independent of $\text{DegNSEP}^*$. 
Thus, steepest descent is generally adaptive to the local strong convexity of the logistic loss function and also achieves faster convergence in a neighborhood of an optimal solution. 
The proofs of our results utilize the generalized self-concordance theory of Bach \cite{bach2010, bach2014adaptivity}, with analysis that is perhaps simpler than the case of stochastic gradient descent examined in \cite{bach2014adaptivity}. 

Theorem \ref{LogitBoost-complexity2} below presents the linear convergence results for steepest descent in the non-separable case. Recall that $\nu^\ast(M)$ denotes the local strong convexity constant of a symmetric positive semidefinite matrix $M$ with respect to the given norm $\|\cdot\|$, and that part {\em (ii)} of Proposition \ref{strictly2} implies that $\nu^\ast(H(\beta^\ast)) > 0$ whenever $\dsep > 0$.

\begin{theorem}\label{LogitBoost-complexity2}{\bf (Linear Convergence of Steepest Descent in the Non-Separable Case)}
Suppose that Steepest Descent (Algorithm \ref{sdm}) is initalized at $\beta^0 := 0$, and is implemented using the step-size sequence:
\begin{equation*}\label{stairmaster2} \alpha_k := \frac{4n\|\nabla L_n(\beta^k)\|_* }{\|\bX\|^2_{\cdot,2}}\ \ \ \mbox{for~all~} k \ge 0 \ .
\end{equation*}
Suppose that the data is non-separable and $\mathrm{DegNSEP}^* > 0$. Define the ``slow" rate of linear convergence constant:
\begin{equation*}
\tau_{\mathrm{slow}} := \left(1 - \frac{2(\mathrm{DegNSEP}^*)\nu^\ast(H(\beta^\ast))n}{(\mathrm{DegNSEP}^* + 2\ln(2)\|\bX\|_{\cdot,\infty})\|\bX\|_{\cdot,2}^2}\right) < 1 \ .
\end{equation*}
Then for all $k \geq 0$, it holds that:
\begin{enumerate}
\item[(i)] (training error): $L_n(\beta^k) - L_n^\ast ~\leq~ (\ln(2) - L_n^\ast)\cdot(\tau_{\mathrm{slow}})^{k}$ , and
\item[(ii)] (coefficient convergence):  $\|\beta^k - \beta^\ast\| ~\leq~ \left(1 + \frac{2\ln(2)\|\bX\|_{\cdot, \infty}}{\mathrm{DegNSEP}^*}\right)\left(\frac{\|\bX\|_{\cdot, 2}}{\nu^\ast(H(\beta^\ast))}\right)\sqrt{\frac{\ln(2) - L_n^\ast}{2n}} \cdot (\tau_{\mathrm{slow}})^{k/2}$ ,
\end{enumerate}
where $\beta^\ast$ is the unique optimal solution of LR.
Furthermore, define:
\begin{equation*}
\check K := \left\lceil\frac{16\ln(2)^2\|\bX\|_{\cdot,2}^4 \|\bX\|_{\cdot,\infty}^2}{9n^2 (\mathrm{DegNSEP}^*)^2\nu^\ast(H(\beta^\ast))^2}\right\rceil \ ,
\end{equation*}
and the ``fast" rate of linear convergence constant:
\begin{equation*}
\tau_{\mathrm{fast}} := \left(1 - \frac{\nu^\ast(H(\beta^\ast))n}{\|\bX\|_{\cdot,2}^2}\right) < \tau_{\mathrm{slow}} < 1\ .
\end{equation*}
Then for all $k \geq \check{K}$, it holds that:
\begin{enumerate}
\item[(iii)] (training error): $L_n(\beta^k) - L_n^\ast ~\leq~ (L_n(\beta^{\check{K}}) - L_n^\ast) \cdot (\tau_{\mathrm{fast}})^{k - \check{K}}$ , and
\item[(iv)] (coefficient convergence):  $\|\beta^k - \beta^\ast\| ~\leq~ \frac{\|\bX\|_{\cdot, 2}}{\nu^\ast(H(\beta^\ast))}\sqrt{\frac{2(L_n(\beta^{\check{K}}) - L_n^\ast)}{n}} \cdot (\tau_{\mathrm{fast}})^{(k - \check{K})/2}$ .
\end{enumerate}
\qed

\end{theorem}\medskip

The proof of Theorem \ref{LogitBoost-complexity2} is presented in Appendix \ref{january3}. As compared to results of a similar flavor for other algorithms, here we have precise guarantees for both the ``slow" and ``fast" rates of linear convergence as well as for the point at which the fast rate is guaranteed to ``kick in." Moreover, there is a natural dependence on $\dsep$ in both the slow rate $\tau_{\mathrm{slow}}$ as well as the iterate $\check K$ by which point the linear convergence rate is improved.
Note that the bound on the training error provided by part {\em (i)} of Theorem \ref{LogitBoost-complexity1}, while sublinear, will be superior to the linear convergence bound provided by part {\em (i)} of Theorem \ref{LogitBoost-complexity2} during the early iterations of steepest descent. On the other hand, the fast linear convergence bound provided by part {\em (iii)} of Theorem \ref{LogitBoost-complexity2} will eventually be the superior of all three bounds when $k$ is large enough.
Recall that Proposition \ref{local_strong_prop} states that $\nu^\ast(H(\beta^\ast))$ is bounded from below as follows:
\begin{equation*}
\nu^\ast(H(\beta^\ast)) ~\geq~ \tfrac{1}{4n}\nu^\ast(\bX^T\bX)\exp \left(-\frac{\ln(2)\|\bX\|_{\cdot,\infty}}{\mathrm{DegNSEP}^*}\right) > 0 \ ,
\end{equation*}
which only depends on $\dsep$ and the magnitude properties of $\bX$. This lower bound can be leveraged to develop versions of the slow and fast linear convergence rates that depend only on $\dsep$ and the magnitude properties of $\bX$. 
However (and as also noted by Bach in \cite{bach2014adaptivity}), the exponential term in the above lower bound is quite pessimistic and in practice $\nu^\ast(H(\beta^\ast))$ tends to be not much smaller than $\tfrac{1}{4n}\nu^\ast(\bX^T\bX)$, i.e, it is as if the exponential term above is not present.\medskip

\subsection{Informing $\ell_2$ Steepest Descent for Solving Logistic Regression in the Separable Case}\label{sdsubsection_sep}
Let us also examine the case when the data is separable.  As mentioned earlier, the logistic regression optimization problem \eqref{poi-logit} is not naturally designed for the case when the data is separable due to the fact that in this case there is no optimal solution.  Indeed, one would suspect that in this case most algorithms -- in particular elementary algorithms such as steepest descent -- would not exhibit computational guarantees of interest.  However, Soudry et al. \cite{soudry2017implicit} and Ji and Telgarsky \cite{telgarsky2018} have recently shown for the case of separable data that steepest descent for the $\ell_2$ norm delivers solutions whose normalized values are approximate-maximum-margin solutions.  The following theorem presents our results for $\ell_2$ steepest descent in the separable case, which adds different and explicit computational guarantees in the separable case.\medskip

\begin{theorem}\label{LogitBoost-complexity3}{\bf ($\ell_2$ Steepest Descent for Logistic Regression:  Separable Case)}
Suppose that $\ell_2$ steepest descent (Algorithm \ref{sdm} with the $\ell_2$ norm) is initialized at $\beta^0 := 0$, and is implemented using the step-size rule:
\begin{equation}\label{stairmaster3} \alpha_k := \frac{2\|\nabla L_n(\beta^k)\|_2}{\|\bX\|_{2,\infty}^2}\ \ \ \mbox{for~all~} k \ge 0 \ .
\end{equation} 

When the data is separable, it holds for all $k \geq 1$ that:

\begin{itemize}

\item[(i)] (margin bound):  there exists $i \in \{0, 1, \ldots, k\}$ for which the normalized iterate $\bar\beta^i := \beta^i/\|\beta^i\|_2$ satisfies 
\begin{equation}\label{new_margin}
\rho(\bar\beta^i) ~\geq~ \frac{\mathrm{DegSEP}^* \cdot \ln\left(\frac{\mathrm{DegSEP}^*}{n\|\bX\|_{2,\infty}}\sqrt{\frac{3(k+1)}{2\ln(2)}} - 1\right)}{2(\ln(k) + 1)} \ ,
\end{equation}

\item[(ii)] (shrinkage): $\|\beta^k\|_2 ~\leq~ \frac{2\ln(k)}{\mathrm{DegSEP}^*} + \frac{2}{\|\bX\|_{2,\infty}}$ , and

\item[(iii)] (gradient bound): $\min\limits_{i \in \{0,\ldots,k\}}\|\nabla L_{n}(\beta^i)\|_2 ~\leq~ \|\bX\|_{2,\infty}\sqrt{\frac{2\ln(2)}{3(k+1)}} \ . $ \ \qed

\end{itemize}

\end{theorem}\medskip

Note that the constant step-size value \eqref{stairmaster3} in Theorem \ref{LogitBoost-complexity3} is different from the corresponding value \eqref{stairmaster1} in Theorem \ref{LogitBoost-complexity1} (when the norm is the $\ell_2$ norm). Indeed, the step-size value \eqref{stairmaster3} is smaller than  the step-size value \eqref{stairmaster1} since $\tfrac{1}{n}\|\bX\|_{2,2}^2 \leq \|\bX\|_{2,\infty}^2$. Therefore, $\tfrac{1}{2}\|\bX\|_{2,\infty}^2$ is a valid upper bound on the Lipschitz constant of the gradient and it is straightforward to develop variants of Theorems \ref{LogitBoost-complexity1} and \ref{LogitBoost-complexity2} for the step-size \eqref{stairmaster3}, wherein $\tfrac{1}{4n}\|\bX\|_{2,2}^2$ would be replaced by $\tfrac{1}{2}\|\bX\|_{2,\infty}^2$ in all of the bounds.

The bound in item {\em (i)} of Theorem \ref{LogitBoost-complexity3} can be understood as $O(1/\ln(k))$ relative convergence to at least $\tfrac{\mathrm{DegSEP}^*}{4}$. To demonstrate this, consider setting $k := \lfloor\frac{2\ln(2)\Omega^2n^2\|\bX\|_{2,\infty}^2}{3(\mathrm{DegSEP}^*)^2}\rfloor$ for some parameter $\Omega \geq 2$. Then the bound in \eqref{new_margin} becomes:
\begin{equation*}
\rho(\bar\beta^i) ~\geq~ \frac{\mathrm{DegSEP}^* \cdot \ln\left(\Omega - 1\right)}{4\ln(\Omega) + 2\ln\left(\frac{2\ln(2)n^2\|\bX\|_{2,\infty}^2}{3(\mathrm{DegSEP}^*)^2}\right) + 2} \ , 
\end{equation*}
and rearranging the above yields:
\begin{align*}
\frac{\rho(\bar\beta^i)}{\tfrac{\mathrm{DegSEP}^*}{4}} ~&\geq~ \frac{\ln\left(\Omega - 1\right)}{\ln(\Omega) + \tfrac{1}{2}\ln\left(\frac{2\ln(2)n^2\|\bX\|_{2,\infty}^2}{3(\mathrm{DegSEP}^*)^2}\right) + \tfrac{1}{2}} \\
~&=~ 1 ~ - ~\frac{\ln\left(\frac{\Omega}{\Omega - 1}\right) + \tfrac{1}{2}\ln\left(\frac{2\ln(2)n^2\|\bX\|_{2,\infty}^2}{3(\mathrm{DegSEP}^*)^2}\right) + \tfrac{1}{2}}{\ln(\Omega) + \tfrac{1}{2}\ln\left(\frac{2\ln(2)n^2\|\bX\|_{2,\infty}^2}{3(\mathrm{DegSEP}^*)^2}\right) + \tfrac{1}{2}} \\
~&\geq~  1 ~ - ~ \frac{\ln(2) + \tfrac{1}{2}\ln\left(\frac{2\ln(2)n^2\|\bX\|_{2,\infty}^2}{3(\mathrm{DegSEP}^*)^2}\right) + \tfrac{1}{2}}{\tfrac{1}{2}\ln(\Omega^2) + \tfrac{1}{2}\ln\left(\frac{2\ln(2)n^2\|\bX\|_{2,\infty}^2}{3(\mathrm{DegSEP}^*)^2}\right) + \tfrac{1}{2}} \\
~&\geq~  1 ~ - ~ \frac{\ln(2) + \tfrac{1}{2}\ln\left(\frac{2\ln(2)n^2\|\bX\|_{2,\infty}^2}{3(\mathrm{DegSEP}^*)^2}\right) + \tfrac{1}{2}}{\tfrac{1}{2}\ln(k+1) + \tfrac{1}{2}} \ ,
\end{align*}
(where the second inequality uses $\Omega \geq 2$), and letting $k$ grow demonstrates $O(1/\ln(k))$ convergence of the iterate margins to $\tfrac{\mathrm{DegSEP}^*}{4}$ or larger.  Interestingly, except for the factor of $4$, this result is similar to Soudry et al. \cite{soudry2017implicit} who show $O(1/\ln(k))$ convergence to $\mathrm{DegSEP}^*$, and to Ji and Telgarsky \cite{telgarsky2018}, whose work shows $O(\sqrt{\ln \ln(k)/\ln(k)})$ convergence to $\mathrm{DegSEP}^*$; however our arguments appear to be entirely different and they result in an explicit margin bound, whereas the results in \cite{telgarsky2018} and \cite{soudry2017implicit} are less transparent.  On the other hand, \cite{telgarsky2018} and \cite{soudry2017implicit} of course prove convergence towards $\mathrm{DegSEP}^*$, not $\frac{\mathrm{DegSEP}^*}{4}$.

In order to prove Theorem \ref{LogitBoost-complexity3}, we will use the following lemma which bounds the norms of iterates of $\ell_2$ steepest descent applied to the logistic regression problem \eqref{poi-logit}, and which is a modified version of a result in Ji and Telgarsky \cite{telgarsky2018}.   

\begin{lemma}{\bf (essentially from Ji and Telgarsky \cite{telgarsky2018})} \label{telgarsky_norm_bound_lemma} Suppose that $\ell_2$ steepest descent (Algorithm \ref{sdm} with the $\ell_2$ norm) is initialized at $\beta^0 := 0$ using the step-size sequence $\{ \alpha_k\}$.  If $\mathrm{DegSEP}^* > 0$ and $\alpha_k \leq \frac{2\|\nabla L_n(\beta^k)\|_2}{\|\bX\|_{2,\infty}^2}$ for all $k \geq 0$, then it holds for all $k \geq 1$ that:
\begin{equation*}
\|\beta^k\|_2 ~\leq~ \frac{2\ln(k)}{\mathrm{DegSEP}^*} + \frac{2}{\|\bX\|_{2,\infty}} \ . 
\end{equation*} \ \qed
\end{lemma}

The proof of this lemma is presented in  Appendix \ref{sweatener}.  

\noindent {\bf Proof of Theorem \ref{LogitBoost-complexity3}:}  
We first prove item {\em (iii)}. Let $i \in \{0, \ldots, k\}$. By the smoothness of the logistic loss function we have:
\begin{align*}
L_n(\beta^{i+1}) ~&\leq~ L_n(\beta^{i}) + \nabla L_n(\beta^i)^T(\beta^{i+1} - \beta^i) + \tfrac{\|\bX\|_{2,2}^2}{8n}\|\beta^{i+1} - \beta^i\|_2^2 \\
~&=~ L_n(\beta^i) - \alpha_i\|\nabla L_n(\beta^i)\|_2^2 + \tfrac{\alpha_i^2\|\bX\|_{2,2}^2}{8n}\|\nabla L_n(\beta^i)\|_2^2 \\
~&=~ L_n(\beta^i) - \tfrac{2}{\|\bX\|_{2,\infty}^2}\|\nabla L_n(\beta^i)\|_2^2 + \tfrac{\|\bX\|_{2,2}^2}{2n\|\bX\|_{2,\infty}^4}\|\nabla L_n(\beta^i)\|_2^2 \\
~&\leq~ L_n(\beta^i) - \tfrac{2}{\|\bX\|_{2,\infty}^2}\|\nabla L_n(\beta^i)\|_2^2 + \tfrac{1}{2\|\bX\|_{2,\infty}^2}\|\nabla L_n(\beta^i)\|_2^2 \\
~&=~ L_n(\beta^i) - \tfrac{3}{2\|\bX\|_{2,\infty}^2}\|\nabla L_n(\beta^i)\|_2^2 \ .
\end{align*}
Summing over all $i \in \{0, \ldots, k\}$ yields:
\begin{equation*}
\frac{3}{2\|\bX\|_{2,\infty}^2}\sum_{i = 0}^k \|\nabla L_n(\beta^i)\|_2^2 \leq L_n(\beta^0) - L_n(\beta^{k+1}) \leq \ln(2) \ ,
\end{equation*}
which implies item {\em (iii)} after rearranging and replacing the terms in the summation by their minimum.  Item {\em (ii)} is a restatement of Lemma \ref{telgarsky_norm_bound_lemma}.  To prove item {\em (i)}, let $i$ be a minimizing index in item {\em (iii)},   and note by item {\em (iii)} and Lemma \ref{margin-lemma} that:
$$
\rho(\beta^i) ~\geq~ \ln\left(\tfrac{\mathrm{DegSEP}^*}{n\|\nabla L_n(\beta^i)\|_2} - 1\right) ~\geq~ \ln\left(\tfrac{\mathrm{DegSEP}^*}{n\|\bX\|_{2,\infty}}\sqrt{\tfrac{3(k+1)}{2\ln(2)}} - 1\right) \ .
$$
Combining the above with item {\em (ii)} and using $\mathrm{DegSEP}^* \leq \|\bX\|_{2,\infty}$, we obtain: 
$$
\rho(\bar\beta^i) = \frac{\rho(\beta^i)}{\|\beta^i\|_2}~\geq~ \frac{\ln\left(\tfrac{\mathrm{DegSEP}^*}{n\|\bX\|_{2,\infty}}\sqrt{\tfrac{3(k+1)}{2\ln(2)}} - 1\right)}{\frac{2\ln(i)}{\mathrm{DegSEP}^*} + \frac{2}{\|\bX\|_{2,\infty}}} ~\ge~\frac{\ln\left(\tfrac{\mathrm{DegSEP}^*}{n\|\bX\|_{2,\infty}}\sqrt{\tfrac{3(k+1)}{2\ln(2)}} - 1\right)}{\frac{2\ln(k)}{\mathrm{DegSEP}^*} + \frac{2}{\mathrm{DegSEP}^*}} \ ,
$$ and the proof then follows from rearranging terms in the above inequality. \qed

While Theorem \ref{LogitBoost-complexity3} holds only for the case of $\ell_2$ steepest descent, it is also possible to derive a much weaker result for any given norm $\|\cdot\|$ by employing the $O(\sqrt{k})$ iterate norm bound in item {\em (ii)} of Theorem \ref{LogitBoost-complexity1} instead of the $O(\ln(k)/\dnsep)$ bound given by Lemma \ref{telgarsky_norm_bound_lemma}. We also mention that Gunasekar et al. \cite{gunasekar2018characterizing} show that $\lim_{k \to \infty} \rho(\bar\beta^k) = \dnsep$ for steepest descent in an arbitrary norm, although there is no analysis of the rate of convergence. Of course, it would be very interesting to see if the tools used to prove these various results can be combined to yield stronger convergence guarantees about the margin.

\section{Informing Stochastic Gradient Descent for Logistic Regression}\label{nye}

In this section, we examine the role of the condition numbers $\dsep$ and $\dnsep$ in the computational and statistical properties of the stochastic gradient descent (SGD) method applied to logistic regression.  Throughout this section, we consider a more general version of the logistic regression problem LR that replaces the empirical average in \eqref{poi-logit} with an arbitrary distribution. Let $\calD$ denote an arbitrary distribution on the data $(\bx, y)$. We consider the following version of the logistic regression problem:
\begin{equation}\label{poi-logit2}
\begin{array}{rccl}
\mathrm{LR}_\calD \ : \ \ \ L_{\calD}^* := & \min\limits_{\beta} &  & L_{\calD}(\beta)  \ := \ \bbE_{(\bx, y) \sim \calD}\left[ \ln\left(1 + \exp\left(-y\beta^T\bx\right) \right) \right] \\ \\
& \text{s.t.} & & \beta \in \mathbb{R}^p \ .
\end{array}
\end{equation}

Note that there are two important special cases of $\mathrm{LR}_\calD$. When $\calD$ is the empirical distribution of the training data $(\bx_1, y_1)\ldots, (\bx_n, y_n)$, then we recover the LR problem \eqref{poi-logit} that has been the primary interest of the previous sections of this paper. On the other hand, it is often useful to conceptualize $\calD$ as the \emph{true} underlying distribution of the data $(\bx, y)$, in which case $\mathrm{LR}_\calD$ is the problem of minimizing the expected logistic loss.  Both of these abstractions will be useful in our analysis of stochastic gradient descent. We refer to $L_{\calD}(\cdot)$ in both cases simply as the ``logistic loss function" with the implicit understanding that the function $L_{\calD}(\cdot)$ is able to capture both the empirical logistic loss as well as the expected logistic loss. Throughout this section, $H(\cdot)$ denotes the Hessian of $L_\calD(\cdot)$.

After briefly reviewing stochastic gradient descent (SGD) for smooth convex optimization in Section \ref{sgd-review}, we examine SGD as applied to logistic regression in two cases: in Section \ref{sect:sgd-nonseparable} we consider the non-separable case and examine $\LRD$ in fully generality, and in Section \ref{sect:sgd-separable} we consider the separable case when $\calD$ is the empirical distribution of a training dataset, i.e., the training data problem LR in \eqref{poi-logit}. 

\subsection{Brief Review of Stochastic Gradient Descent}\label{sgd-review}
We first review the stochastic gradient descent method for solving the generic unconstrained differentiable convex optimization problem \eqref{poi4}.  By way of motivating context, it is sometimes the case that computing the gradient $\nabla f(\cdot)$ of $f(\cdot)$ is very expensive or even intractable, but it may be relatively easy to compute a stochastic estimate of $\nabla f(x)$ at $x$, which we denote by $\tg f(x)$, via a stochastic gradient oracle.  We say that the oracle computes an unbiased stochastic gradient if $\bbE[\tg f(x) ~|~ x] = \nabla f(x)$.  Notice that by construction $\tg f(x)$ is a conditional random variable given $x$. 
The basic stochastic gradient descent (SGD) method is presented in Algorithm \ref{sgd}, whose structure is the same as steepest descent (Algorithm \ref{sdm}) for the $\ell_2$ norm, the only difference being that the stochastic gradient $\tg f(x^i)$ replaces the exact gradient $\nabla f(x^i)$ in Step (1.).  Also, for simplicity, we only consider the case of a constant step-size sequence $\alpha_i := \bar\alpha > 0$ for all $i \geq 0$. 

\floatname{algorithm}{Algorithm}
\begin{algorithm}
\caption{Stochastic Gradient Descent with constant step-size $\bar\alpha$}\label{sgd}
\begin{algorithmic}
\STATE {\bf Initialize} at $x^0 \in \mathbb{R}^p$, $i \leftarrow 0$, and set the total number of iterations $k \geq 1$. \medskip

{\bf At iteration $i$:}
\STATE 1. Call stochastic oracle to compute $\tg f(x^i)$
\STATE 2. Update:
\begin{description}
\item $\ \  x^{i+1} \gets x^i - \bar\alpha \cdot \tg f(x^i)$
\end{description}
\medskip

{\bf After $k$ iterations:}
\STATE Option A:  Output $\hat{x}^k \gets \frac{1}{k+1}\sum_{i = 0}^k x^i$.
\STATE Option B:  Output $\hat{x}^k \gets x^{I_k}$ where $I_k$ is a random variable distributed uniformly on $\{0, 1, \ldots, k\}$.
\end{algorithmic}
\end{algorithm}

Stochastic gradient descent, and more generally the idea of stochastic approximation, dates back to the seminal work of Robbins and Monro \cite{robbins1951}. For recent works related to stochastic gradient descent see, e.g., Nemirovski et al. \cite{nemirovski2009robust}, Bottou \cite{bottou2010large}, Lan \cite{lan2012optimal}, Bottou et al. \cite{bottou2018optimization}, and the references therein.
The output of Algorithm \ref{sgd} using either Option A or Option B depends on the entire sequence $x^0, x^1, \ldots, x^k$. Nevertheless in both options $\hat x^k$ can be updated in an online fashion which does not require storage of the entire sequence $x^0, x^1, \ldots, x^k$. Clearly in the case of Option A we have $\hat{x}_k = \frac{k}{k+1}\hat{x}^{k-1} + \frac{1}{k+1}x^k$ for $k \geq 1$, where by convention $\hat{x}^{0} = x^0$. In the case of Option B, note that we can construct $I_k$ recursively as follows: given $I_{k-1}$, which is uniformly distributed on $\{0, \ldots, k - 1\}$, we define $I_k$ to be equal to $I_{k-1}$ with probability $\frac{k}{k+1}$ and equal to $k$ with probability $\frac{1}{k+1}$. Then it holds that $I_k$ is uniformly distributed on $\{0, 1, \ldots, k\}$ and is independent of $x^0, x^1, \ldots, x^k$. As pointed out in \cite{ghadimi2013stochastic}, Option B can also be implemented by first generating $I_k$ uniformly at random on $\{0, \ldots, k\}$ during the initialization stage of the algorithm, and then only running for $I_k$ iterations before stopping (assuming $k$ is fixed in advance).

In addition to the smoothness condition \eqref{smoothness} (with respect to the $\ell_2$ norm in this case), a condition that is required in the typical analysis of SGD is the following bounded second moment condition:
\begin{equation}\label{bounded_variance}
\bbE\left[\|\tg f(x)\|_2^2 ~|~ x\right] \leq M^2  \ \ \ \text{for all} \ x\in \mathbb{R}^p \ ,
\end{equation}
where $M$ is a positive constant. In the case of smooth convex optimization, as is studied in \cite{lan2012optimal} for example, \eqref{bounded_variance} is often replaced with a bound on the variance of the stochastic gradient oracle instead, which is smaller. However for our purposes in studying logistic regression, \eqref{bounded_variance} is adequate and simplifies the analysis.

The following theorem presents a stylized version of computational guarantees for SGD, which is an amended and extended version of mostly well-known results about stochastic gradient descent, see for example Nemirovski et al. \cite{nemirovski2009robust} as well as Ghadimi and Lan \cite{ghadimi2013stochastic}. In the theorem, the expectation in items {\em (i)} and {\em (iii)} is taken with respect to all of the stochasticity of Algorithm \ref{sgd}. \medskip

\begin{theorem}\label{sgdcomplexity}{\bf (Computational Guarantees for Stochastic Gradient Descent)} Let $\{x^k\}$ be generated according to the stochastic gradient descent method (Algorithm \ref{sgd}) using a constant step-size $\bar\alpha > 0$.  Under either Option A or Option B, it holds for all $k \geq 0$ that:

\begin{itemize}
\item[(i)] (expected optimality gap): $$\bbE[f(\hat{x}^k)] - f(x) ~\leq~ \displaystyle\frac{\|x^0 - x\|_2^2}{2\bar\alpha (k+1)} + \frac{\bar\alpha M^2}{2} \ \  \text{ for all } x \in \bbR^p \  , \ \mbox{and}$$
\item[(ii)] (norm bound):  $$\|x^k - x^0\|_2 ~\leq~ \bar\alpha\sum_{i = 0}^{k-1}\|\tg f(x^i)\|_2 \ .$$ 
\end{itemize}
Under Option B, it holds for all $k \geq 0$ that:
\begin{itemize}
\item[(iii)] (expected gradient bound):  $$\bbE\left[\|\nabla f(\hat{x}^k)\|_2^2\right] ~\leq~ \frac{f(x^0) - f^\ast}{\bar\alpha(k+1)} + \frac{\bar\alpha L M^2}{2} \ .$$
\end{itemize} \qed
\end{theorem}

For completeness, a self-contained proof of Theorem \ref{sgdcomplexity} is given in Appendix \ref{unitedflight}. Note that when an optimal solution $x^\ast$ of \eqref{poi4} exists, then we can take $x \gets x^\ast$ in item {\em (i)} to obtain a bound on the expected optimality gap $\bbE[f(\hat{x}^k)] - f^\ast$. Items {\em (i)} and {\em (iii)} of Theorem \ref{sgdcomplexity} present bounds that hold in expectation; bounds that hold with high probability require additional assumptions such as moment generating function type assumptions, compactness of the feasible region, etc., see \cite{nemirovski2009robust} and \cite{ghadimi2013stochastic}.

\subsection{Informing Stochastic Gradient Descent for Logistic Regression in the Non-Separable Case}\label{sect:sgd-nonseparable}
Let us now return to the logistic regression problem $\LRD$, which we examine in full generality in the case of non-separable data. First we need to extend the definitions of non-separable and separable datasets to an arbitrary distribution $\calD$ over the data $(\bx, y)$.
Let $\supp(\calD) \subseteq \bbR^p \times \{-1, +1\}$ denote the support of the distribution $\calD$. Then we say that the data distribution $\calD$ is \emph{separable} if there exists a model $\beta \in \bbR^p$ such that $\inf_{(\bx, y) \in \supp(\calD)}y \beta^T\bx > 0$. Otherwise, if $\inf_{(\bx, y) \in \supp(\calD)}y \beta^T\bx \leq 0$ for every model $\beta$, then we say that the data distribution $\calD$ is \emph{non-separable}.

As in the previously examined case of finite training datasets, clearly some non-separable distributions might be ``more non-separable'' than others, so let us introduce a way to measure the extent to which the distribution is non-separable.  We define the degree of non-separability of the distribution $\calD$ (with respect to the norm $\|\cdot\|$) to be:
\begin{equation}\label{dsepD}
\begin{array}{rccl}
 \ \ \ \mathrm{DegNSEP}^*_\calD := & \displaystyle\min_{\beta \in \mathbb{R}^p} &  & \bbE_{(\bx, y) \sim \calD}\left[[y\beta^T\bx]^{-}\right] \\ \\ 
& \text{s.t.} & & \|\beta\| = 1 \ ,
\end{array}
\end{equation}
which states that $\dsepD$ is the smallest (over all normalized models $\beta$) expected misclassification error of the model $\beta$.  Note that the norm $\|\cdot\|$ in the above definition is any generic given norm. It is straightforward to extend Proposition \ref{strictly2} to this more general setting, and we present this generalization in Proposition \ref{strictlyD}. 
\begin{proposition}\label{strictlyD}  If $\dsepD > 0$, then:
\begin{enumerate}
\item[(i)] there is a unique optimal solution $\beta^\ast$ of the logistic regression problem $\LRD$,
\item[(ii)] $H(\beta^\ast) \succ 0$ ,
\item[(iii)] $\|\beta^*\| ~\le~ \displaystyle\frac{L_\calD^\ast}{\dsepD} ~\le~ \displaystyle\frac{\ln(2)}{\dsepD}$ , and
\item[(iv)] $\mathrm{Dist}_0 ~\le~ \displaystyle\frac{\ln(2) + L_\calD^\ast}{\dsepD} ~\le~ \displaystyle\frac{2\ln(2)}{\dsepD}$ .
\end{enumerate}
\end{proposition}
\noindent {\bf Proof:} The proof follows the same structure as the proof of Proposition \ref{strictly2}, but requires a more careful argument to prove that $H(\beta) \succ 0$ for all $\beta \in \bbR^p$. Suppose that $\dsepD > 0$, and let $\beta \in \bbR^p$ be given. It can also be demonstrated (e.g., by Section 7.2.4 of \cite{shapiro2009lectures}) that $L_\calD(\cdot)$ is twice differentiable and that $H(\beta) = \bbE_{(\bx, y) \sim \calD}\left[\ell^{\prime\prime}(y\beta^T\bx)\bx\bx^T\right]$. Note that $H(\beta) \succeq 0$ by convexity. Suppose, by way of contradiction, that there exists $\bar\beta \in \bbR^p$ with $\|\bar\beta\| = 1$ and $\bar\beta^TH(\beta)\bar\beta = 0$. Then we have that $0 = \bar\beta^TH(\beta)\bar\beta = \bbE_{(\bx, y) \sim \calD}\left[\ell^{\prime\prime}(y\beta^T\bx)\bar\beta^T\bx\bx^T\bar\beta\right] = \bbE_{(\bx, y) \sim \calD}\left[\ell^{\prime\prime}(y\beta^T\bx)(\bar\beta^T\bx)^2\right]$. Therefore it holds that $\ell^{\prime\prime}(y\beta^T\bx)(\bar\beta^T\bx)^2 = 0$ with probability one, and since $\ell^{\prime\prime}(y\beta^T\bx) > 0$ we must have that $\bar\beta^T\bx = 0$ with probability one. This then implies that $\bbE_{(\bx, y) \sim \calD}\left[[y\bar\beta^T\bx]^{-}\right] = 0$, which implies that $\dsepD = 0$, and this provides the desired contradiction. Therefore $H(\beta) \succ 0$. The remainder of the proof exactly follows that of Proposition \ref{strictly2}, and is omitted for brevity.\qed

Throughout this section we assume the following finiteness condition on second moments of $\calD$:

\begin{assumption}{\bf (Finite second moments of the data distribution $\calD$)}\label{Dassumption}
The data distribution $\calD$ satisfies $\bbE\left[\|\bx\|_2^2\right] < +\infty$.\end{assumption}

Define the second moment matrix $\Sigma := \bbE[\bx\bx^T] \in \bbR^{p \times p}$.  It follows from Assumption \ref{Dassumption} that $\Sigma$ is well-defined and finite, and that $L_\calD(\cdot)$ is continuous, convex, differentiable, and satisfies $\nabla L_\calD(\beta) = \bbE\left[\nabla_\beta \ln\left(1 + \exp\left(-y\beta^T\bx\right) \right)\right]$ for all $\beta \in \bbR^p$ (see, e.g., \cite{shapiro2009lectures}).   We also have:

\begin{proposition}{\bf{(Lipschitz smoothness of the logistic loss function)}}\label{LcalD}
The logistic loss function $L_{\calD}(\cdot)$ is $L = \tfrac{1}{4}\lambda_{\max}(\Sigma)$-smooth with respect to the $\ell_2$ norm on $\bbR^p$. \ \qed
\end{proposition}

A proof of Proposition \ref{LcalD} is given in Appendix \ref{subway}.

Denote the scalar logistic loss function by $\ell(t) := \ln(1 + \exp(-t))$.  We assume the following regarding the stochastic gradient oracle for the logistic loss function $L_\calD(\beta)$ of \eqref{poi-logit2}:

\begin{assumption}{\bf (Stochastic gradient oracle for the logistic loss function)}\label{oracleAssumption}
The stochastic gradient oracle $\tg L_\calD(\cdot)$ is implemented by drawing an independent sample from the distribution $\calD$.  That is, for any (possibly random) $\beta \in \bbR^p$, the stochastic gradient $\tg L_\calD(\beta)$ is computed by independently sampling $(\tilde\bx, \tilde y)$ from the distribution $\calD$ and assigning $\tg L_\calD(\beta) \gets \nabla_\beta \ell(\tilde y \beta^T\tilde\bx)$.
\end{assumption}

\begin{proposition}{\bf (Second moment of the stochastic gradient)}\label{logit-var}
The stochastic gradient $\tg L_\calD(\cdot)$ computed from the oracle described in Assumption \ref{oracleAssumption} satisfies the second moment upper bound \eqref{bounded_variance} with $M^2 = \tr(\Sigma)$. \ \qed
\end{proposition}

A proof of Proposition \ref{logit-var} is given in Appendix \ref{valencia_cafesoret}.  

Theorem \ref{LogitSGDComplexity1} presents the main computational guarantees for SGD applied to the logistic regression problem \eqref{poi-logit2} in the case when the data distribution $\calD$ is non-separable.

\begin{theorem}\label{LogitSGDComplexity1}{\bf (Computational Guarantees for SGD:  Non-Separable Case)}
Suppose that SGD (Algorithm \ref{sgd}) for logistic regression is initialized at $\beta^0 := 0$ and is implemented using the constant step-size $ \bar\alpha > 0$,
and that the stochastic gradients are computed as in Assumption \ref{oracleAssumption}.
Under either Option A or Option B of SGD (Algorithm \ref{sgd}), if the data distribution $\calD$ is non-separable, then for all $k \geq 0$ it holds that:
\begin{equation}\label{sgd_upper_bound}
\bbE\left[L_\calD(\hat\beta^k)\right] - L_\calD^\ast ~\leq~ \frac{(\ln(2))^2}{2\bar\alpha (k + 1) \cdot (\dsepD)^2} ~+~ \frac{\bar\alpha \cdot \tr(\Sigma)}{2} \ .
\end{equation}
\end{theorem}

\noindent {\bf Proof:}
This result is a straightforward application of Theorem \ref{sgdcomplexity}, Proposition \ref{strictlyD}, and Proposition \ref{logit-var}. By item {\em (i)} of Proposition \ref{strictlyD} an optimal solution $\beta^\ast$ exists. The result follows by directly applying item {\em (i)} of Theorem \ref{sgdcomplexity}, along with the bound $\|\beta^*\| \le \frac{\ln(2)}{\dsepD}$ provided by item {\em (ii)} of Proposition \ref{strictlyD}, and the value of $M^2 = \tr(\Sigma)$ from Proposition \ref{logit-var}.
\qed

Theorem \ref{LogitSGDComplexity1} provides an upper bound on the expected logistic loss that naturally depends on the condition number $\dsepD$ and that holds for any step-size value $\bar\alpha$. Given knowledge of the constants $\dsepD$ and $\tr(\Sigma)$, it is possible to tune $\bar\alpha$ in order to minimize the upper bound expression on the right side of \eqref{sgd_upper_bound} for a given $k$. However, in practice one does not typically know either of these constants.  Indeed, it is more realistic to assume one has knowledge of a deterministic upper bound on the size of the feature vectors, i.e., there is an available constant $R > 0$ for which $\|\bx\|_2 \leq R$ with probability one.  Under this assumption it also follows that $\tr(\Sigma) = \bbE[\|\bx\|_2^2] \leq R^2$.  Corollary \ref{Rcor} presents a computational guarantee for SGD in the case of non-separable data under a slightly weaker assumption and using a step-size that only incorporates knowledge of the constant $R$.

\begin{corollary}\label{Rcor}
Suppose that we have available a constant $R$ such that $\tr(\Sigma) \leq R^2$.
Consider running SGD (Algorithm \ref{sgd}) for a total of $k$ iterations, initialized at $\beta^0 := 0$, with stochastic gradients computed as in Assumption \ref{oracleAssumption}, and using the constant step-size 
\begin{equation*}
 \bar{\alpha} := \frac{\ln(2)}{R^2\sqrt{k+1}} \ . 
 \end{equation*}
Under either Option A or Option B of SGD (Algorithm \ref{sgd}), if the data distribution $\calD$ is non-separable it holds for all $k \ge 0$ that:
\begin{equation*}
\bbE\left[L_\calD(\hat\beta^k)\right] - L_\calD^\ast ~\leq~ \frac{\ln(2)}{2\sqrt{k+1}}\left(\frac{R^2}{(\dsepD)^2} + 1\right) \ .
\end{equation*} \qed
\end{corollary}

Theorem \ref{LogitSGDComplexity1} and Corollary \ref{Rcor} present results that highlight the role of the condition number $\dsepD$ in the well-known $O(1/\sqrt{k})$ computational guarantee for stochastic gradient descent. It is also possible to study how $\dsepD$ informs the adaptive $O(1/\mu k)$ (where $\mu$ is the smallest eigenvalue of the Hessian at the optimum) guarantees developed by Bach in \cite{bach2014adaptivity}. Proposition \ref{bach_for_sgd} below directly follows from Proposition 10 of \cite{bach2014adaptivity} and demonstrates how $\dsepD$ informs the corresponding computational guarantees developed therein. Note that $\lambda_{\min}(H(\beta^\ast)) > 0$ by part {\em (ii)} of Proposition \ref{strictlyD}. Note that the step-size considered in Proposition \ref{bach_for_sgd} is larger than that considered in \cite{bach2014adaptivity} by a constant factor of $2\ln(2)$. 
\begin{proposition}\label{bach_for_sgd}{\bf (Bach \cite{bach2014adaptivity}, Proposition 10)}
Suppose that we have available a constant $R$ such that $\|\bx\|_2 \leq R$ with probability one.
Consider running SGD (Algorithm \ref{sgd}) for a total of $k$ iterations, initialized at $\beta^0 := 0$, with stochastic gradients computed as in Assumption \ref{oracleAssumption}, and using the constant step-size 
\begin{equation*}
 \bar{\alpha} := \frac{\ln(2)}{R^2\sqrt{k+1}} \ . 
\end{equation*}
Under Option A of SGD (Algorithm \ref{sgd}), if the data distribution $\calD$ is non-separable and $\mathrm{DegNSEP}^*_\calD > 0$, then it holds for all $k \ge 0$ that:
\begin{enumerate}
\item[(i)] \ \ \ $\displaystyle\bbE\left[L_\calD(\hat\beta^k)\right] - L_\calD^\ast ~\leq~ \frac{R^2}{\lambda_{\min}(H(\beta^\ast))(k+1)}\left(\frac{10R(\ln(2))^2}{\dsepD} + 15\right)^4$ , and
\item[(ii)] \ \ \ $\displaystyle\bbE\left[\|\hat\beta^k - \beta^\ast\|_2^2\right] ~\leq~ \frac{R^2}{\lambda_{\min}(H(\beta^\ast))^2(k+1)}\left(\frac{12R(\ln(2))^2}{\dsepD} + 21\right)^4$ ,
\end{enumerate}
where $\beta^\ast$ is the unique optimal solution of $\mathrm{LR}_\calD$.
\end{proposition}
\noindent {\bf Proof:}
Consider the scaled objective function $\tilde{L}_\calD(\cdot) := 2\ln(2)L_\calD(\cdot)$, and denote its Hessian by $\tilde H(\cdot)$. It can be verified that $\tilde{L}_\calD(\cdot)$ satisfies assumptions (A1) - (A7) of \cite{bach2014adaptivity} with constant $2\ln(2)R$. 
Clearly, running SGD with objective $L_\calD(\cdot)$ and constant step-size $\frac{\ln(2)}{R^2\sqrt{k+1}}$ is equivalent to running SGD with objective $\tilde{L}_\calD(\cdot)$ and constant step-size $\frac{1}{2R^2\sqrt{k+1}}$, which is required by Proposition 10 of \cite{bach2014adaptivity}. Therefore, directly applying this proposition to $\tilde{L}_n(\cdot)$ yields:
\begin{equation*}
\bbE\left[\tilde{L}_\calD(\hat\beta^k)\right] - \tilde{L}_\calD^\ast ~\leq~ \frac{(2\ln(2))^2R^2}{\lambda_{\min}(\tilde{H}(\beta^\ast))(k+1)}\left(10R\ln(2)\|\beta^\ast\|_2 + 15\right)^4 \ .
\end{equation*} 
Finally, using $\tilde{H}(\beta^\ast) = 2\ln(2) H(\beta^\ast)$, the upper bound $\|\beta^*\| \le \frac{\ln(2)}{\dsepD}$ provided by item {\em (ii)} of Proposition \ref{strictlyD}, and dividing both sides of the above by $2\ln(2)$ yields part {\em (i)} of the theorem. Part {\em (ii)} also directly follows from Proposition 10 of \cite{bach2014adaptivity} by a similar argument.
\qed

Comparing Corollary \ref{Rcor} with Proposition \ref{bach_for_sgd}, note that the relative sizes of the constants $R$, $\dsepD$, and $\lambda_{\min}(H(\beta^\ast))$, as well as the total number of iterations $k$ will determine which bound dominates the other. And of course if $k$ is large enough then Proposition \ref{bach_for_sgd} yields the better bound.

\subsection{Informing Stochastic Gradient Descent for Logistic Regression in the Separable Case}\label{sect:sgd-separable}
In this subsection we examine the properties of SGD in the case when the data is separable, for the previously examined logistic regression problem \eqref{poi-logit}, i.e., we assume that $\calD$ is the empirical distribution of the training data $(\bx_1, y_1)\ldots, (\bx_n, y_n)$ in this subsection and we revert to the notation used throughout Section \ref{imp}.
The following theorem presents our results in this case.

\begin{theorem}\label{sgd_separable}{\bf (Computational Guarantees for SGD:  Separable Case)}
Suppose that we have available a constant $R$ such that $\|\bx_i\|_2 \leq R$ for all $i \in \{1, \ldots, n\}$.
Consider running SGD (Algorithm \ref{sgd}) for a total of $k$ iterations, initialized at $\beta^0 := 0$, with stochastic gradients computed as in Assumption \ref{oracleAssumption}, and using the constant step-size 
\begin{equation*}
\bar{\alpha} := \frac{\ln(2)}{R^2\sqrt{k+1}} \ . 
\end{equation*}
Under Option B of SGD (Algorithm \ref{sgd}), when the data is separable it holds for all $k \ge 1$ that:
\begin{itemize}

\item[(i)] (margin bound): For any $\gamma \in (0,1]$, with probability at least $1 - \gamma$ the normalized iterate $\bar{\beta}^k := \hat\beta^k/\|\hat\beta^k\|$ satisfies
\begin{equation}\label{sgd_margin_bound}
\rho(\bar\beta^k) ~>~ \frac{\mathrm{DegSEP}^* \cdot \ln\left(\frac{\mathrm{DegSEP}^*\sqrt{\gamma}\sqrt[4]{k+1}}{nR\sqrt{1.1}} - 1\right)}{2(\ln(k) + 1)}
\end{equation}

\item[(ii)] (shrinkage): $\displaystyle\|\hat\beta^k\|_2 ~\leq~ \frac{2\ln(k)}{\mathrm{DegSEP}^*} + \frac{2}{\|\bX\|_{2,\infty}} \ , $ and

\item[(iii)] (expected gradient bound): $\displaystyle\bbE\left[\|\nabla L_{n}(\hat\beta^k)\|_2^2\right] ~<~ \frac{1.1 \cdot R^2}{\sqrt{k+1}} \ .$  \qed

\end{itemize}
\end{theorem}

The margin bound in Theorem \ref{sgd_separable} is similar in flavor to the bound for steepest descent in Theorem \ref{LogitBoost-complexity3}, but is weaker (due to stochasticity). Indeed, by similar arguments as in Section \ref{sdsubsection_sep}, the bound in item {\em (i)} of Theorem \ref{sgd_separable} implies a computational guarantee of the form
\begin{equation}\label{eps_guarantee}
\frac{\rho(\bar\beta^k)}{\tfrac{\mathrm{DegSEP}^*}{8}} \geq 1 - \frac{C}{\ln(k+1)} \ \text{ with probability at least $1 - \gamma$} \ ,
\end{equation}
for any fixed $\gamma \in (0,1]$, with $C = 4\ln(2) - 2 \ln(\gamma) + \ln(\tfrac{(1.1)^2n^4R^4}{(\mathrm{DegSEP}^*)^4}) + 1$. (This should be compared to the case of deterministic steepest descent where we have deterministic convergence to at least $\tfrac{\mathrm{DegSEP}^*}{4}$.) To demonstrate this, consider setting $k := \lfloor\frac{\Omega^4 (1.1)^2n^4R^4}{(\mathrm{DegSEP}^*)^4\gamma^2}\rfloor$ for some parameter $\Omega \geq 2$. Then the bound in \eqref{sgd_margin_bound} becomes:
\begin{equation*}
\rho(\bar\beta^k) ~\geq~ \frac{\mathrm{DegSEP}^* \cdot \ln\left(\Omega - 1\right)}{8\ln(\Omega) - 4\ln(\gamma) + 2\ln\left(\frac{(1.1)^2n^4R^4}{(\mathrm{DegSEP}^*)^4}\right) + 2} \ , 
\end{equation*}
and rearranging the above and using $\Omega \geq 2$ yields:
\begin{align*}
\frac{\rho(\bar\beta^k)}{\tfrac{\mathrm{DegSEP}^*}{8}}
~&\geq~  1 ~ - ~ \frac{\ln(2) - \tfrac{1}{2}\ln(\gamma) + \tfrac{1}{4}\ln\left(\frac{(1.1)^2n^4R^4}{(\mathrm{DegSEP}^*)^4}\right) + \tfrac{1}{4}}{\ln(\Omega) - \tfrac{1}{2}\ln(\gamma) + \tfrac{1}{4}\ln\left(\frac{(1.1)^2n^4R^4}{(\mathrm{DegSEP}^*)^4}\right) + \tfrac{1}{4}} \\
~&\geq~  1 ~ - ~ \frac{\ln(2) - \tfrac{1}{2}\ln(\gamma) + \tfrac{1}{4}\ln\left(\frac{(1.1)^2n^4R^4}{(\mathrm{DegSEP}^*)^4}\right) + \tfrac{1}{4}}{\tfrac{1}{4}\ln(k+1)} ~ = ~ 1 - \frac{C}{\ln(k+1)}  \ .
\end{align*}This result is similar to results in \cite{srebro2018b}, wherein $O(1/\ln(k))$ convergence towards $\mathrm{DegSEP}^*$ is demonstrated for SGD for sampling without replacement.  Note that we provide an explicit margin bound in item {\em (i)} above and that we study SGD with sampling with replacement.
On the other hand, \cite{srebro2018b} of course proves convergence towards $\mathrm{DegSEP}^*$, not $\frac{\mathrm{DegSEP}^*}{8}$.  It would be interesting to see if the tools used to prove all of these results can be combined somehow to yield stronger convergence guarantees about the margin for SGD. Note also that the margin bound \eqref{sgd_margin_bound} is proven by applying Markov's inequality with the bound on the second moment of $\|\nabla L_{n}(\hat\beta^k)\|_2$ given by item {\em (iii)} of the theorem. In the case of non-separable data, Bach \cite{bach2014adaptivity} is able to strengthen this second moment bound to $O(1/k)$ and also derives bounds on the higher-order moments of $\|\nabla L_{n}(\hat\beta^k)\|_2$ (for Option A of SGD). It would also be interesting to see if similar bounds can be derived and used to strengthen the margin bound in the case of separable data.

In order to prove Theorem \ref{sgd_separable}, we will use the following lemma, which is similar to Lemma \ref{telgarsky_norm_bound_lemma} and bounds the norms of iterates of SGD applied to the logistic regression problem \eqref{poi-logit}.

\begin{lemma}{\bf (essentially from Ji and Telgarsky \cite{telgarsky2018})} \label{telgarsky_norm_bound_lemma_sgd} Suppose that SGD (Algorithm \ref{sgd}) is initialized at $\beta^0 := 0$ using the constant step-size value $\bar\alpha$.  If $\mathrm{DegSEP}^* > 0$ and $\bar\alpha \leq \frac{2}{\|\bX\|_{2,\infty}^2}$, then it holds for all $k \geq 1$ that:
\begin{equation*}
\|\beta^k\|_2 ~\leq~ \frac{2\ln(k)}{\mathrm{DegSEP}^*} + \frac{2}{\|\bX\|_{2,\infty}} \ . 
\end{equation*} \ \qed
\end{lemma}

The proof of this lemma is presented in  Appendix \ref{sweatener}.  

\noindent {\bf Proof of Theorem \ref{sgd_separable}:}
Item {\em (ii)} follows directly from Lemma \ref{telgarsky_norm_bound_lemma_sgd} as well as the fact that $\ln(\cdot)$ is an increasing function.  Item {\em (iii)} is a straightforward application of item {\em (iii)} of Theorem \ref{sgdcomplexity}. Indeed, recalling that $L_n(\beta^0) - L_n^\ast = \ln(2)$ (equality holds in the separable case) and $\lambda_{\max}(\Sigma) \leq \tr(\Sigma) \leq R^2$, item {\em (iii)} follows directly from the definition of $\bar\alpha$, Propositions \ref{LcalD} and \ref{logit-var}, and item {\em (iii)} of Theorem \ref{sgdcomplexity}.

To prove item {\em (i)} first note that Markov's inequality yields:
\begin{equation*}
\bbP\left(\|\nabla L_{n}(\hat\beta^k)\|_2^2 \geq \frac{1.1 \cdot R^2}{\gamma\sqrt{k+1}}\right) \leq \frac{\bbE\left[\|\nabla L_{n}(\hat\beta^k)\|_2^2\right]}{\frac{1.1 \cdot R^2}{\gamma\sqrt{k+1}}} < \gamma \ ,
\end{equation*}
where the second inequality follows from item {\em (iii)} of the theorem. Therefore with probability at least $1 - \gamma$ it holds that:
\begin{equation}\label{grad_markov}
\|\nabla L_{n}(\hat\beta^k)\|_2^2 < \frac{1.1 \cdot R^2}{\gamma\sqrt{k+1}} \ .
\end{equation}
We will now demonstrate that if \eqref{grad_markov} holds (in addition to everything else) then \eqref{sgd_margin_bound} holds, which therefore implies that the statement in part {\em (i)} is true. Indeed, combining \eqref{grad_markov} with Lemma \ref{margin-lemma} yields:
\begin{equation*}
\rho(\hat\beta^k) ~\geq~ \ln\left(\tfrac{\mathrm{DegSEP}^*}{n\|\nabla L_n(\hat\beta^k)\|_2} - 1\right) ~ > ~ \ln\left(\tfrac{\mathrm{DegSEP}^*\sqrt{\gamma}\sqrt[4]{k+1}}{nR\sqrt{1.1}} - 1\right) \ .
\end{equation*}
Combining the above with item {\em (ii)} and using $\mathrm{DegSEP}^* \leq \|\bX\|_{2,\infty}$, we obtain:
\begin{equation*}
\rho(\bar\beta^k) = \frac{\rho(\hat\beta^k)}{\|\hat\beta^k\|_2} ~ > ~ \frac{\ln\left(\tfrac{\mathrm{DegSEP}^*\sqrt{\gamma}\sqrt[4]{k+1}}{nR\sqrt{1.1}} - 1\right)}{\frac{2\ln(k)}{\mathrm{DegSEP}^*} + \frac{2}{\|\bX\|_{2,\infty}}} ~\ge~ \frac{\ln\left(\tfrac{\mathrm{DegSEP}^*\sqrt{\gamma}\sqrt[4]{k+1}}{nR\sqrt{1.1}} - 1\right)}{\frac{2\ln(k)}{\mathrm{DegSEP}^*} + \frac{2}{\mathrm{DegSEP}^*}} \ ,
\end{equation*}
and the proof then follows from rearranging terms in the above inequality.
\qed

\section{Conclusions}
The theme of this paper is the interplay between data conditioning, behavior/properties of the optimization problem, and computational guarantees of first-order methods, all in the context of logistic regression. We have presented results that make rigorous the intuitive notion that the optimization problem itself as well as the corresponding algorithms for training a logistic regression model are well-behaved when the degree of non-separability of the dataset is large. We also have presented results that demonstrate that the specific algorithmic properties of steepest descent and stochastic gradient descent lead to large margin solutions in the case of separable data, which runs counter to the intuition that logistic regression is ill-behaved in this case. We hope that further examination of the role of data and problem conditioning in the analysis of other statistical learning problems and other algorithms will extend the general understanding of these problems and algorithms.

\bibliographystyle{abbrv}
\bibliography{GF-papers-orc_student_paper}

\appendix

\section{Appendix}

\subsection{Proofs of Propositions \ref{jimrenegar1} and \ref{jimrenegar2}}\label{June29}

\noindent {\bf Proof of Proposition \ref{jimrenegar1}:}  Notice that the optimization problem defining \eqref{dsep} has a continuous objective function and a compact feasible region, whereby it follows from the Weierstrass Theorem that \eqref{dsep} attains its optimum at some $\bar\beta$ and therefore $\dsep = \tfrac{1}{n}\sum_{i=1}^n [y_i \bx_i^T\bar\beta]^-$.  It follows from norm duality that there exists $\bar s$ satisfying $\|\bar s\|_* = 1$ and $\bar s^T\bar \beta = \|\bar\beta\| = 1$.  Let $\varepsilon >0$ be given, and now define $\Delta \bX := u\bar s^T$ where $u_i := y_i[y_i \bx_i^T\bar\beta]^- + y_i \varepsilon $.  Notice for each $i=1, \ldots, n$ that $y_i(\bx_i + \Delta \bx_i)^T\bar\beta = y_i \bx_i^T\bar\beta + y_i u_i \bar s^T\bar \beta = y_i \bx_i^T\bar\beta + [y_i \bx_i^T\bar\beta]^- + \varepsilon \ge \varepsilon >0$, whereby the perturbed dataset $(\bX + \Delta \bX, y)$ is separable and hence $\pertsep \le \tfrac{1}{n}\|\Delta \bX\|_{\cdot, 1} = \tfrac{1}{n}\|u\|_1 \|\bar s\|_* = \tfrac{1}{n}\|u\|_1 =  \dsep + \varepsilon$.  As this is true for any $\varepsilon >0$ it follows that  $\pertsep \le \dsep $.

We next show that $\dsep \le \pertsep$, which will complete the proof.  Suppose that $\Delta \bX$ satisfies $(\bX + \Delta \bX, y)$ is separable, and hence there exists $\beta$ with $\|\beta\| = 1$ and $y_i(\bx_i + \Delta \bx_i)^T\beta >0$ for $i=1, \ldots, n$.  Define the vector $v$ component-wise for $i=1, \ldots, n$ by:
$$v_i := \left\{ \begin{array}{ll} 0 \ \ \ & \mbox{if}  \ \ y_i\bx_i^T\beta \ge 0 \\
y_i  & \mbox{if}  \ \ y_i\bx_i^T\beta < 0  \ , \end{array} \right. $$  
and notice in particular that if $y_i\bx_i^T\beta < 0$ then $[y_i\bx_i^T\beta]^- = -y_i\bx_i^T\beta < y_i(\Delta \bx_i)^T\beta = v_i(\Delta \bx_i)^T\beta$.  Also, if $y_i\bx_i^T\beta \ge 0$, then $[y_i\bx_i^T\beta]^- = 0 = v_i(\Delta \bx_i)^T\beta$.  Therefore $\dsep \le \tfrac{1}{n}\sum_{i=1}^n [y_i\bx_i^T\beta]^- \le  \tfrac{1}{n}\sum_{i=1}^n v_i(\Delta \bx_i)^T\beta = \tfrac{1}{n} v^T\Delta \bX \beta \le \tfrac{1}{n}\|\Delta \bX\|_{\cdot, 1}$ since $\|v\|_\infty \le 1$.  Thus $\dsep \le \tfrac{1}{n}\|\Delta \bX\|_{\cdot,1}$ for any perturbation $\Delta \bX$ for which $(\bX + \Delta \bX, y)$ is separable, and hence $\dsep \le \pertsep$, completing the proof. \qed

\noindent {\bf Proof of Proposition \ref{jimrenegar2}:}  Define $\Delta_n = \{ \lambda \in \mathbb{R}^n : e^T\lambda = 1, \ \lambda \ge 0 \}$. We can write \eqref{dnsep} in maxmin form as:
\begin{equation}\label{cold}
\dnsep :=  \displaystyle\max_{\beta: \|\beta\| \le 1}   \  \displaystyle\min_{\lambda \in \Delta_n} \lambda^T Y \bX \beta \ \ = \ \  \displaystyle\min_{\lambda \in \Delta_n} \ \displaystyle\max_{\beta: \|\beta\| \le 1}   \ \lambda^T Y \bX \beta \ \ = \ \   \displaystyle\min_{\lambda \in \Delta_n} \|\bX^T Y \lambda\|_* \ , 
\end{equation}
where the middle equality follows from minmax strong duality.  Furthermore, both the minmax problem and the maxmin problem attain their optima for some $\bar\beta$ satisfying $\|\bar\beta\| \le 1$ and $\bar\lambda \in \Delta_n$ which implies that:
\begin{equation}\label{icy}
\dnsep \ := \ \bar\lambda^T Y \bX \bar\beta  \ = \   \rho(\bar\beta) \ = \min_{i} (y_i \bx_i^T \bar\beta)_i \ = \ \|\bX^T Y \bar\lambda\|_* \ . 
\end{equation}
Now define $\Delta \bX := -y\bar\lambda^TY\bX$.  Direct substitution yields $\bar\lambda^TY(\bX+\Delta \bX)=0$, which then implies that there does not exist any $\beta$ satisfying $Y(\bX+ \Delta \bX)\beta >0$, and hence $(\bX + \Delta \bX, y)$ is not separable.  Therefore $\pertnsep \le \|\Delta \bX\|_{\cdot, \infty} = \| -y\bar\lambda^TY\bX \|_{\cdot , \infty} = \|y\|_\infty\|\bX^TY\bar\lambda\|_* = \dnsep$.

We next show that $\dnsep \le \pertnsep$, which will complete the proof.  Suppose that $\Delta \bX$ satisfies $(\bX + \Delta \bX, y)$ is not separable, and hence by a theorem of the alternative there exists $\lambda \in \Delta_n$ satisfying $\lambda^TY(\bX + \Delta \bX) = 0 $.  Using the values $\bar\beta$ and $\bar\lambda$ defined above, we have:

$$\dnsep = \bar\lambda^TY\bX\bar\beta \le \lambda^TY\bX\bar\beta = -\lambda^TY\Delta\bX\bar\beta \le \|\Delta \bX\|_{\cdot, \infty} \|\bar\beta\| \|Y\lambda\|_1 \le \|\Delta \bX\|_{\cdot, \infty}    \ . $$

Thus $\dnsep \le \|\Delta \bX\|_{\cdot,\infty}$ for any perturbation $\Delta \bX$ for which $(\bX + \Delta \bX, y)$ is non-separable, and hence $\dnsep \le \pertnsep$, completing the proof. \qed

\subsection{Proof of Theorem \ref{sdmcomplexity}}\label{january2}

Since $f(\cdot)$ satisfies \eqref{smoothness}, it follows easily from the fundamental theorem of calculus that:
\begin{equation}\label{smooth-ineq}
f(y) \leq f(x) + \nabla f(x)^T(y - x) + \tfrac{L}{2}\|y-x\|^2 \ \ \ \mbox{for~all~} x,y \ .
\end{equation}  (For a short proof of this fact, see Proposition A.2 of \cite{fg2013} for example.)
Applying \eqref{smooth-ineq} to the iterates of the Steepest Descent Method yields the following for each $i \geq 0$:
\begin{equation}\label{running}\begin{array}{rcl}
f(x^{i+1}) & \leq &f(x^{i}) + \nabla f(x^{i})^T(x^{i+1} - x^{i}) + \tfrac{L}{2}\|x^{i+1} - x^{i}\|^2\\ \\
&=& f(x^{i}) - \alpha_i \|\nabla f(x^i)\|_*+ \tfrac{L}{2}\alpha_i^2 \ ,
\end{array}\end{equation}
where the equality follows since $\|\nabla f(x^k)\|_\ast = \max\limits_{d : \|d\| \leq 1}\nabla f(x^k)^Td = \nabla f(x^k)^T d^k $.  Summing the above for $i=0, \ldots, k$ yields:
\begin{equation}\label{running2}\begin{array}{rcl}f^* \le f(x^{k+1}) \le f(x^0) -\sum_{i=0}^k \alpha_i \|\nabla f(x^i)\|_* + \tfrac{L}{2}\sum_{i=0}^k \alpha_i^2 \ . \end{array}\end{equation}
Next notice that $$\sum_{i=0}^k \alpha_i \|\nabla f(x^i)\|_* \ge \left(\sum_{i=0}^k \alpha_i \right) \left(\min_{i \in \{0,\ldots,k\}} \|\nabla f(x^i)\|_*\right) \ , $$ and substituting this inequality above and rearranging yields \begin{equation}\label{sdmarbitrary}
\min\limits_{i \in \{0,\ldots,k\}}\|\nabla f(x^i)\|_* \leq  \frac{f(x^0) - f^* +\frac{L}{2}\sum_{i = 0}^{k}\alpha_i^2}{\sum_{i=0}^{k}\alpha_i} \ .
\end{equation}
Now suppose we use the step-sizes \eqref{sdmquadrule}.  Substituting \eqref{sdmquadrule} into \eqref{running} yields:
\begin{equation}\label{descent-ineq}
f(x^{i+1}) \leq f(x^{i}) - \tfrac{1}{2L}\|\nabla f(x^i)\|_*^2 \ ,
\end{equation}
which shows that the values $f(x^i)$ are monotone decreasing and hence $f(x^i) \le f(x^0)$, whereby $x^i \in {\cal S}_0$.  Substituting the step-sizes \eqref{sdmquadrule} into \eqref{running2} yields after rearranging:
\begin{equation}\label{wb}
\sum_{i = 0}^k \|\nabla f(x^i)\|_*^2 \ \leq 2L(f(x^0) - f(x^{k+1})) \leq \ 2L(f(x^0) - f^* ) \ ,
\end{equation}
and therefore
\begin{equation*}
(k+1)\left(\min_{i \in \{0,\ldots,k\}}\|\nabla f(x^i)\|_*\right)^2  \leq \sum_{i = 0}^k \|\nabla f(x^i)\|_*^2 \ \leq \ 2L(f(x^0) - f^* ) \ ,
\end{equation*}
and rearranging yields {\em (iv)}.  Now suppose as well that $\mathrm{Dist}_0$ is finite, and let $x^i$ be an iterate of the steepest descent method.  It was shown above that $x^i \in {\cal S}_0$,  whereby there exists $x^* \in {\cal S}^*$ for which $\|x^i -x^*\| \le \mathrm{Dist}_0$, and from the gradient inequality for the convex function $f(\cdot)$ it holds that $$\begin{array}{rcl} f^* = f(x^*) &\ge& f(x^i) + \nabla f(x^i)^T(x^*-x^i) \\ \\
&\ge&  f(x^i) -\| \nabla f(x^i)\|_* \|x^*-x^i\| \\ \\
&\ge&  f(x^i) -\| \nabla f(x^i)\|_* \mathrm{Dist}_0 \ ,
\end{array}$$ and rearranging the above yields $\| \nabla f(x^i)\|_* \ge \frac{f(x^i)-f^*}{ \mathrm{Dist}_0}$.  Substituting this inequality into \eqref{descent-ineq} and subtracting $f^*$ from both sides yields:
$$f(x^{i+1}) -f^*\leq f(x^{i}) -f^*- \frac{(f(x^{i}) -f^*)^2}{2L\mathrm{Dist}_0^2} \ . $$
Define $a_i := f(x^{i}) -f^*$, and it follows that the nonnegative series $\{a_i\}$ satisfies $a_{i+1} \le a_i - \frac{a_i^2}{2L\mathrm{Dist}_0^2}$.  A standard induction argument (see for example Lemma 3.5 of \cite{beckgcd}) then establishes that $$a_{k} \le \frac{1}{\frac{1}{a_0} + \frac{k}{2L\mathrm{Dist}_0^2}} \ , $$ which when rearranged yields the first inequality of {\em (i)}.  The second inequality of {\em (i)} follows since $\hat K^0 > 0$.

Rearranging \eqref{descent-ineq} yields $\|\nabla f(x^i)\|_*^2 \le 2L(f(x^i) - f(x^{i+1}) \le 2L(f(x^i) - f^*)$, which after taking square roots proves the first inequality of {\em (ii)}, and the second inequality of {\em (ii)} follows by substituting in the bound on $f(x^k) - f^*$ from the first inequality of {\em (i)}.

To prove {\em (iii)}, use $k-1$ in \eqref{wb}, and use the step-lengths \eqref{sdmquadrule} to yield:
$$2L(f(x^0)-f^*) \ge \sum_{i = 0}^{k-1} \|\nabla f(x^i)\|_*^2 = L^2 \sum_{i = 0}^{k-1} \alpha_i^2 \ge \left(\frac{1}{k}\right) L^2  \left(\sum_{i = 0}^{k-1} \alpha_i \right)^2 \ge \left(\frac{1}{k}\right) L^2 \|x^k - x^0\|^2 \ , $$
and rearranging the above yields {\em (iii)}. \qed\medskip

\subsection{Proof of Proposition \ref{logistic-props}}\label{tired}

We first present a property of the following ``prox'' function $d(\cdot): [0,1]^n \rightarrow \mathbb{R}$ defined by:

\begin{equation}\label{dee} d(w) :=  \frac{1}{n}\left[\sum_{i=1}^n w_i \ln(w_i) + (1-w_i) \ln (1-w_i) \right] \ ,  \end{equation}
where $\alpha \ln(\alpha) := 0$ for $\alpha = 0$. \medskip

\begin{proposition}\label{appendix1} Consider the function $d(\cdot): [0,1]^n \rightarrow \mathbb{R}$ given by \eqref{dee}.  It holds that $d(\cdot)$ is a $\sigma := \frac{4}{n}$-strongly convex function with respect to the Euclidean norm $\|w\|:=\|w\|_2$.
\end{proposition}

\noindent {\bf Proof:}  Let $G:=[0,1]^n$, consider any point $w \in \tint G$, and let $H(w)$ denote the Hessian matrix of $d(\cdot)$ at $w$.  The off-diagonal components of $H(w)$ are all zero, and the $i^{\mathrm{th}}$ diagonal component is $H_{ii}(w) = \frac{1}{n}\frac{1}{w_i(1-w_i)} \ge \frac{4}{n}$, and hence $v^T[H(w)]v \ge \frac{4}{n}v^Tv$ for any $v$.  Now let $y \in \tint G$.  Invoking an intermediate value theorem of calculus, there exists a scalar $c \in [0,1]$ for which it holds that:
$$ d(y) = d(w) + \nabla d(w)^T(y-w) + \tfrac{1}{2}(y-w)H(w+c(y-w))(y-w) \ , $$
whereby:
$$ d(y) \ge d(w) + \nabla d(w)^T(y-w) + \tfrac{1}{2}(y-w)\left[\tfrac{4}{n}I\right](y-w) = d(w) + \nabla d(w)^T(y-w) + \tfrac{1}{2}\tfrac{4}{n}\|y-w\|_2^2 \ . $$  This proves that $d(\cdot)$ is $\sigma=\frac{4}{n}$-strongly convex on $\tint G$, and a continuity argument establishes the result for all of $G$. \qed

\noindent {\bf Proof of Proposition \ref{logistic-props}:}  We first claim that
\begin{equation}\label{logistic-rep}
L_{n}(\beta) = \max\limits_{w \in [0,1]^n}\left\{-w^T[\tfrac{1}{n}Y\bX]\beta - d(w)\right\} \end{equation}
where $d(\cdot)$ is given by \eqref{dee}, and the unique optimal solution to the maximization problem in \eqref{logistic-rep} is:
\begin{equation}\label{logistic-sol}
w^\ast(\beta)_i = \frac{1}{1+\exp(y_i \beta^T\bx_i)} \ \ ,  \  i = 1, \ldots, n \ .
\end{equation}Indeed, it is easy to verify through optimality conditions that the unique optimal solution to the maximization problem in \eqref{logistic-rep} is given by \eqref{logistic-sol}, and direct substitution and simplification of terms then yields the equality in \eqref{logistic-rep}.  Using the representation \eqref{logistic-rep}, Theorem 1 of ~\cite{nest05smoothing} implies that the Lipschitz constant $L$ of the gradient of $L_{n}(\cdot)$ is at most $[\frac{1}{n}\|\bX\|_{\cdot,2}]^2/\sigma$ where $\sigma$ is the strong convexity parameter of the function $d(\cdot)$.  From Proposition \ref{appendix1} it holds that $\sigma \ge \frac{4}{n}$, which implies that $L \le [\frac{1}{n}\|\bX\|_{\cdot,2}]^2/\sigma \le [\frac{1}{4n}]\|\bX\|_{\cdot,2}^2$. \qed

\subsection{Proof of Theorem \ref{LogitBoost-complexity2}}\label{january3}
Following Bach \cite{bach2010, bach2014adaptivity}, a three-times differentiable convex function $f(\cdot) : \mathbb{R}^p \to \mathbb{R}$ is said to be generalized self-concordant (with respect to the norm $\|\cdot\|$) if there is a constant $R > 0$ such that for all $x, \hat x \in \bbR^p$, the scalar function $\varphi(\cdot) : t \mapsto f(x + t(\hat x - x))$ satisfies:
\begin{equation}\label{gen_self_con_f}
|\varphi^{\prime\prime\prime}(t)| ~\leq~ R\|x - \hat x\| \varphi^{\prime\prime}(t) \ \text{ for all } t \in \bbR \ .
\end{equation}
The above definition is a slight modification of that given in \cite{bach2010, bach2014adaptivity}, which works with the $\ell_2$ norm. The following proposition, which is a very minor generalization of a result shown in \cite{bach2010}, demonstrates that the logistic loss function $L_n(\cdot)$ is generalized self-concordant with constant $R = \|\bX\|_{\cdot, \infty}$.

\begin{proposition}\label{logistic_gen_sc}
The logistic loss function $L_n(\cdot)$ is generalized self-concordant (with respect to the norm $\|\cdot\|$) with constant $R = \|\bX\|_{\cdot, \infty}$.
\end{proposition}
\noindent {\bf Proof:}
Let $\beta, \hat \beta \in \bbR^p$ be given and define the scalar function $\varphi(\cdot) : \bbR \to \bbR$ by $\varphi(t) := L_n(\beta + t(\hat \beta - \beta))$. A simple calculation yields:
\begin{align*}
|\varphi^{\prime\prime\prime}(t)| ~&=~ \left|\frac{1}{n}\sum_{i = 1}^n \ell^{\prime\prime\prime}(y_i\beta^T\bx_i + ty_i(\hat{\beta} - \beta)^T\bx_i)\cdot(y_i(\hat{\beta} - \beta)^T\bx_i)^3\right| \\
~&\leq~ \frac{1}{n}\sum_{i = 1}^n |\ell^{\prime\prime\prime}(y_i\beta^T\bx_i + ty_i(\hat{\beta} - \beta)^T\bx_i)|\cdot(y_i(\hat{\beta} - \beta)^T\bx_i)^2 \cdot |y_i(\hat{\beta} - \beta)^T\bx_i| \\
~&\leq~ \frac{1}{n}\sum_{i = 1}^n \ell^{\prime\prime}(y_i\beta^T\bx_i + ty_i(\hat{\beta} - \beta)^T\bx_i)\cdot(y_i(\hat{\beta} - \beta)^T\bx_i)^2 \cdot |y_i(\hat{\beta} - \beta)^T\bx_i| \\
~&\leq~ \frac{1}{n}\sum_{i = 1}^n \ell^{\prime\prime}(y_i\beta^T\bx_i + ty_i(\hat{\beta} - \beta)^T\bx_i)\cdot(y_i(\hat{\beta} - \beta)^T\bx_i)^2 \cdot \|\bx_i\|_\ast \|\beta - \hat\beta\| \\
~&\leq~ \frac{\|\bX\|_{\cdot, \infty}\|\beta - \hat\beta\|}{n}\sum_{i = 1}^n \ell^{\prime\prime}(y_i\beta^T\bx_i + ty_i(\hat{\beta} - \beta)^T\bx_i)\cdot(y_i(\hat{\beta} - \beta)^T\bx_i)^2 \\
~&=~ \|\bX\|_{\cdot, \infty}\|\beta - \hat\beta\| \varphi^{\prime\prime}(t) \ ,
\end{align*}
where the second inequality above uses $|\ell^{\prime\prime\prime}(\cdot)| \leq \ell^{\prime\prime}(\cdot)$, the third uses H{\"o}lder's inequality, and the final inequality uses $\|\bX\|_{\cdot, \infty} = \max\limits_{i \in \{1, \ldots, n\}} \|\bx_i\|_\ast$.
\qed

In order to prove Theorem \ref{LogitBoost-complexity2}, we use the following lemma which is a minor extension of Lemma 9 of \cite{bach2014adaptivity}.
\begin{lemma}\label{bach_lemma9}{\bf (essentially Bach \cite{bach2014adaptivity}, Lemma 9)}
Suppose that $\mathrm{DegNSEP}^* > 0$. Let $\beta$ satisfying $L_n(\beta) \le \ln(2)$ be given, and let $\beta^*$ be the unique optimal solution of LR. Then it holds that:
\begin{enumerate}
\item[(i)] \ \ \ $L_n(\beta) - L_n^\ast ~\leq~ \left(1 + \frac{2\ln(2)\|\bX\|_{\cdot, \infty}}{\mathrm{DegNSEP}^*}\right)\frac{\|\nabla L_n(\beta)\|_\ast^2}{\nu^\ast(H(\beta^\ast))}$ , and
\item[(ii)] \ \ \ $\|\beta - \beta^\ast\| ~\leq~ \left(1 + \frac{2\ln(2)\|\bX\|_{\cdot, \infty}}{\mathrm{DegNSEP}^*}\right)\left(\frac{\|\bX\|_{\cdot, 2}}{\nu^\ast(H(\beta^\ast))}\right)\sqrt{\frac{L_n(\beta) - L_n^\ast}{2n}}$ .
\end{enumerate}
If in addition $\beta$ satisfies $\frac{\|\nabla L_n(\beta)\|_\ast \|\bX\|_{\cdot, \infty}}{\nu^\ast(H(\beta^\ast))} ~\leq~ \frac{3}{4}$, then it holds that:
\begin{enumerate}
\item[(iii)] \ \ \ $L_n(\beta) - L_n^\ast ~\leq~ \frac{2\|\nabla L_n(\beta)\|_\ast^2}{\nu^\ast(H(\beta^\ast))}$ , and
\item[(iv)] \ \ \ $\|\beta - \beta^\ast\| ~\leq~ \frac{\|\bX\|_{\cdot, 2}}{\nu^\ast(H(\beta^\ast))}\sqrt{\frac{2(L_n(\beta) - L_n^\ast)}{n}}$ \ .
\end{enumerate}
\end{lemma}
\noindent {\bf Proof:}
First note that if $\beta = \beta^\ast$ then the lemma is trivial, so we assume that $\beta \neq \beta^\ast$. 
We will apply Lemma 13 of \cite{bach2014adaptivity} to the scalar function $\varphi(\cdot) : [0,1] \to \bbR$ defined by $\varphi(t) := L_n(\beta^\ast + t(\beta - \beta^\ast))$. 
Let us define $S := \|\bX\|_{\cdot, \infty}\|\beta - \beta^\ast\|$, and note that Proposition \ref{logistic_gen_sc} implies that $\varphi(\cdot)$ satisfies $|\varphi^{\prime\prime\prime}(t)| \leq S\varphi^{\prime\prime}(t)$ for all $t \in [0,1]$.
Simple calculations yield:
\begin{equation*}
\varphi^\prime(t) = \nabla L_n(\beta^\ast + t(\beta - \beta^\ast))^T(\beta - \beta^\ast) \ \text{ and } \ \varphi^{\prime\prime}(t) = (\beta - \beta^\ast)^T H(\beta^\ast + t(\beta - \beta^\ast))(\beta - \beta^\ast) \ \text{ for all } t \in [0,1] \ .
\end{equation*}
In particular, we have $\varphi^{\prime}(0) = 0$ by the optimality of $\beta^\ast$ and
\begin{equation}\label{varphi_ineq1}
\varphi^\prime(1) = \nabla L_n(\beta)^T(\beta - \beta^\ast) \leq \|\nabla L_n(\beta)\|_\ast \|\beta - \beta^\ast\| \ ,
\end{equation}
by H{\"o}lder's inequality. Moreover, we have:
\begin{equation}\label{varphi_ineq2}
\varphi^{\prime\prime}(0) = (\beta - \beta^\ast)^T H(\beta^\ast)(\beta - \beta^\ast) \geq \nu^\ast(H(\beta^\ast))\|\beta - \beta^\ast\|^2 > 0 \ ,
\end{equation}
by the definition of $\nu^\ast(H(\beta^\ast))$ in \eqref{local_strong_constant}, part {\em (ii)} of Proposition \ref{strictly2}, and since $\beta \neq \beta^\ast$.
Therefore $\varphi(\cdot)$ satisfies the hypotheses of Lemma 13 of \cite{bach2014adaptivity}, and a direct application of this lemma yields:
\begin{equation}\label{bach_conclusion}
\frac{\varphi^{\prime}(1)}{\varphi^{\prime\prime}(0)}S \geq 1 - \exp(-S) \ , \text{ and } \ \varphi(1) \leq \varphi(0) + \frac{\varphi^{\prime}(1)^2}{\varphi^{\prime\prime}(0)}(1 + S) \ .
\end{equation}
Following the proof of Proposition \ref{strictly2}, we have that $S = \|\bX\|_{\cdot, \infty}\|\beta - \beta^\ast\| \leq \frac{2\ln(2)\|\bX\|_{\cdot, \infty}}{\mathrm{DegNSEP}^*}$ since $\beta$ satisfies $L_n(\beta) \leq \ln(2)$. Substituting this upper bound on $S$ along with the inequalities \eqref{varphi_ineq1} and \eqref{varphi_ineq2} into the rightmost inequality in \eqref{bach_conclusion} yields part {\em (i)} of the lemma. Making the same substitutions into the leftmost inequality in \eqref{bach_conclusion} and rearranging yields:
\begin{equation*}
\|\beta - \beta^\ast\| \leq \frac{S\|\nabla L_n(\beta)\|_\ast}{(1 - \exp(-S))\nu^\ast(H(\beta^\ast))} \leq \left(1 + \frac{2\ln(2)\|\bX\|_{\cdot, \infty}}{\mathrm{DegNSEP}^*}\right)\frac{\|\nabla L_n(\beta)\|_\ast}{\nu^\ast(H(\beta^\ast))} \ ,
\end{equation*}
where the second inequality uses $\frac{S}{1 - \exp(-S)} \leq 1 + S$. Applying \eqref{descent-ineq} along with Proposition \ref{logistic-props} in this context yields $L_n^\ast \leq L_n(\beta) - \tfrac{2n}{\|\bX\|_{\cdot, 2}^2}\|\nabla L_n(\beta)\|_\ast^2$, 
which after rearranging terms and combining with the above inequality yields part {\em (ii)} of the lemma.

Now suppose that $\frac{\|\nabla L_n(\beta)\|_\ast \|\bX\|_{\cdot, \infty}}{\nu^\ast(H(\beta^\ast))} \leq \frac{3}{4}$ additionally holds.  Then: 
\begin{equation*}
\frac{\varphi^{\prime}(1)S}{\varphi^{\prime\prime}(0)} \leq \frac{\|\nabla L_n(\beta)\|_\ast \|\beta - \beta^\ast\|^2 \|\bX\|_{\cdot, \infty}}{\nu^\ast(H(\beta^\ast))\|\beta - \beta^\ast\|^2} = \frac{\|\nabla L_n(\beta)\|_\ast \|\bX\|_{\cdot, \infty}}{\nu^\ast(H(\beta^\ast))} \leq \frac{3}{4} \ .
\end{equation*}
Hence following Lemma 13 of \cite{bach2014adaptivity}, it holds that $\varphi^{\prime\prime}(0) \leq 2\varphi^{\prime}(1)$ and $\varphi(1) \leq \varphi(0) + 2\frac{\varphi^{\prime}(1)^2}{\varphi^{\prime\prime}(0)}$, and making the same substitutions as before yields parts {\em (iii)} and {\em (iv)}. 
\qed

\noindent {\bf Proof of Theorem \ref{LogitBoost-complexity2}:}
Let $k \geq 0$ be given. Applying \eqref{descent-ineq} along with Proposition \ref{logistic-props} in this context yields:
\begin{equation}\label{descent_logistic}
L_n(\beta^{k+1}) \leq L_n(\beta^k) - \tfrac{2n}{\|\bX\|_{\cdot,2}^2}\|\nabla L_n(\beta^k)\|_\ast^2 \ .
\end{equation}
In particular, $L_n(\beta^{k+1}) \leq L_n(\beta^k)$, which implies that $L_n(\beta^i) \leq L_n(\beta^0) = \ln(2)$ for all $i \geq 0$.
Therefore, we may apply item {\em (i)} of Lemma \ref{bach_lemma9}, which yields:
\begin{equation*}
L_n(\beta^k) - L_n^\ast ~\leq~ \left(1 + \frac{2\ln(2)\|\bX\|_{\cdot, \infty}}{\mathrm{DegNSEP}^*}\right)\frac{\|\nabla L_n(\beta^k)\|_\ast^2}{\nu^\ast(H(\beta^\ast))} \ .
\end{equation*}
Combining the above inequality with \eqref{descent_logistic} yields:
\begin{equation*}
L_n(\beta^{k+1}) \leq L_n(\beta^k) - \frac{2n\nu^\ast(H(\beta^\ast))(L_n(\beta^k) - L_n^\ast)}{\|\bX\|_{\cdot,2}^2\left(1 + \frac{2\ln(2)\|\bX\|_{\cdot, \infty}}{\mathrm{DegNSEP}^*}\right)} \ . 
\end{equation*}
Finally, subtracting $L_n^\ast$ from both sides of the above and rearranging terms yields:
\begin{equation*}
L_n(\beta^{k+1}) - L_n^\ast \leq (L_n(\beta^k) - L_n^\ast)\left(1 - \frac{2(\mathrm{DegNSEP}^*)\nu^\ast(H(\beta^\ast))n}{(\mathrm{DegNSEP}^* + 2\ln(2)\|\bX\|_{\cdot, \infty})\|\bX\|_{\cdot,2}^2}\right) \ ,
\end{equation*}
which immediately implies part {\em (i)} of the theorem. Part {\em (ii)} of the theorem follows by substituting the bound on $L_n(\beta^k) - L_n^\ast$ from part {\em (i)} of the theorem into the bound on $\|\beta^k - \beta^\ast\|$ from part {\em (ii)} of Lemma \ref{bach_lemma9}. 

Now assume that $k \geq \check K$. Part {\em (iii)} of Theorem \ref{LogitBoost-complexity1} implies that:
\begin{equation*}
\frac{\|\nabla L_n(\beta^{k})\|_\ast \|\bX\|_{\cdot, \infty}}{\nu^\ast(H(\beta^\ast))} ~\leq~ \frac{\|\bX\|_{\cdot,2}^2\ln(2)\|\bX\|_{\cdot, \infty}}{\nu^\ast(H(\beta^\ast))\sqrt{k} \cdot n \cdot \mathrm{DegNSEP}^*} \leq \frac{3}{4} \ , 
\end{equation*}
where the final inequality uses the fact that $k \geq \check K \geq \frac{16\ln(2)^2\|\bX\|_{\cdot,2}^4 \|\bX\|_{\cdot,\infty}^2}{9n^2 (\mathrm{DegNSEP}^*)^2\nu^\ast(H(\beta^\ast))^2}$. 
Thus we may apply part {\em (iii)} of Lemma \ref{bach_lemma9}, which yields $L_n(\beta^k) - L_n^\ast \leq \tfrac{2\|\nabla L_n(\beta^k)\|_\ast^2}{\nu^\ast(H(\beta^\ast))}$.
By the same arguments as above, we obtain:
\begin{equation*}
L_n(\beta^{k+1}) - L_n^\ast \leq (L_n(\beta^k) - L_n^\ast)\left(1 - \frac{\nu^\ast(H(\beta^\ast))n}{\|\bX\|_{\cdot,2}^2}\right) \ ,
\end{equation*}
which immediately implies part {\em (iii)} of the thoerem. Part {\em (iv)} of the theorem similarly follows by substituting the bound on $L_n(\beta^k) - L_n^\ast$ from part {\em (iii)} into the bound on $\|\beta^k - \beta^\ast\|$ from part {\em (iv)} of Lemma \ref{bach_lemma9}.
\qed

\subsection{Proofs of Lemmas \ref{telgarsky_norm_bound_lemma} and \ref{telgarsky_norm_bound_lemma_sgd}}\label{sweatener}

Denote the univariate logistic loss function by $\ell(t):= \ln(1+ \exp(-t))$.  We start with the following quite general proposition which presents a bound on the iterate sequence $\{\beta_k\}$ of \emph{any} algorithm whose step direction $g_k$ is an average of gradients of the logistic loss function over a subset $S_k$ of the observations, for all $k$.

\begin{proposition}\label{telgarsky_prop} Consider any algorithm for solving the logistic regression problem \eqref{poi-logit}, and let $\{\beta^k\}$ denote the iterate sequence.  Suppose that $\beta^0:=0$ and $\{\beta^k\}$ satisfies:
\begin{equation}\label{iterates}
\beta^{k+1} = \beta^k - \alpha_k g_k \ \text{ where } \ g_k = \frac{1}{|S_k|}\sum_{i \in S_k} \nabla_\beta \ell(y_i(\beta^k)^T\bx_i) \ 
\end{equation}
for some $S_k \subseteq \{1, \ldots, n\}$, for all $k \ge 0$. If $\alpha_j \leq \frac{2}{\|\bX\|_{2,\infty}^2}$ for all $j \geq 0$, then for all $\beta \in \bbR^p$ and for all $k \geq 0$ it holds that:
\begin{equation*}
\|\beta^k - \beta\|_2^2 ~\leq~ \|\beta\|_2^2 + 2\sum_{j = 0}^{k - 1}\frac{\alpha_j}{|S_j|}\sum_{i \in S_j}\ell(y_i\beta^T\bx_i) \ .
\end{equation*}
\end{proposition}
\noindent {\bf Proof:}
For any $j \in \{0, \ldots, k-1\}$ it holds that:
\begin{align*}
\|\beta^{j+1} - \beta\|_2^2 &= \|\beta^{j} - \beta - \alpha_j g_j\|_2^2 \\
&= \|\beta^{j} - \beta\|_2^2 - 2\alpha_j g_j^T(\beta^j - \beta) + \alpha_j^2\|g_j\|_2^2 \\
&\leq \|\beta^{j} - \beta\|_2^2 + 2\tfrac{\alpha_j}{|S_j|}\textstyle\sum_{i \in S_j} [\ell(y_i\beta^T\bx_i) - \ell(y_i(\beta^j)^T\bx_i)] + \alpha_j^2\|g_j\|_2^2 \\
&= \|\beta^{j} - \beta\|_2^2 + 2\tfrac{\alpha_j}{|S_j|}\textstyle\sum_{i \in S_j} [\ell(y_i\beta^T\bx_i) - \ell(y_i(\beta^j)^T\bx_i)] + \tfrac{\alpha_j^2}{|S_j|^2}\|\textstyle\sum_{i \in S_j} \nabla_\beta \ell(y_i(\beta^j)^T\bx_i)\|_2^2 \\
&\leq \|\beta^{j} - \beta\|_2^2 + 2\tfrac{\alpha_j}{|S_j|}\textstyle\sum_{i \in S_j} [\ell(y_i\beta^T\bx_i) - \ell(y_i(\beta^j)^T\bx_i)] + \tfrac{\alpha_j^2}{|S_j|^2}\left(\textstyle\sum_{i \in S_j}\|\nabla_\beta \ell(y_i(\beta^j)^T\bx_i)\|_2\right)^2 \\
&\leq \|\beta^{j} - \beta\|_2^2 + 2\tfrac{\alpha_j}{|S_j|}\textstyle\sum_{i \in S_j} [\ell(y_i\beta^T\bx_i) - \ell(y_i(\beta^j)^T\bx_i)] + \tfrac{\alpha_j^2}{|S_j|}\textstyle\sum_{i \in S_j}\|\nabla_\beta \ell(y_i(\beta^j)^T\bx_i)\|_2^2 \\
&= \|\beta^{j} - \beta\|_2^2 + 2\tfrac{\alpha_j}{|S_j|}\textstyle\sum_{i \in S_j}\ell(y_i\beta^T\bx_i) + \tfrac{\alpha_j}{|S_j|}\textstyle\sum_{i \in S_j}[\alpha_j\|\nabla_\beta \ell(y_i(\beta^j)^T\bx_i)\|_2^2 - 2\ell(y_i(\beta^j)^T\bx_i)] \ ,
\end{align*}
where the first inequality above is an application of the gradient inequality, the second inequality above uses the triangle inequality, and the third inequality utilizes an inequality between the $\ell_1$ and $\ell_2$ norms. Now since $\alpha_j \leq \frac{2}{\|\bX\|_{2,\infty}^2}$, it holds that:
\begin{equation*}
\alpha_j\|\nabla_\beta \ell(y_i(\beta^j)^T\bx_i)\|_2^2 = \alpha_j\ell^\prime(y_i(\beta^j)^T\bx_i)^2\|\bx_i\|_2^2 \leq 2\ell^\prime(y_i(\beta^j)^T\bx_i)^2 \leq 2\ell(y_i(\beta^j)^T\bx_i)) \ ,
\end{equation*}
where the last inequality uses $\ell^{\prime}(\cdot)^2 \leq |\ell^\prime(\cdot)| \leq \ell(\cdot)$. Therefore:
\begin{equation*}
\|\beta^{j+1} - \beta\|_2^2 ~\leq~ \|\beta^{j} - \beta\|_2^2 + 2\frac{\alpha_j}{|S_j|}\sum_{i \in S_j}\ell(y_i\beta^T\bx_i) \ ,
\end{equation*}
and summing the previous inequality over $j \in \{0, \ldots, k-1\}$ yields the result.\qed

\noindent {\bf Proofs of Lemmas \ref{telgarsky_norm_bound_lemma} and \ref{telgarsky_norm_bound_lemma_sgd}:}
We present a unified proof of these two results. Let $\bar{\beta}$ denote the normalized maximum margin hyperplane, i.e., the optimal solution of \eqref{dnsep}, and define $\tilde{\beta}^k := (\ln(k)/\mathrm{DegSEP}^*)\bar\beta$. For each $i \in \{1, \ldots, n\}$, we have:
\begin{equation*}
y_i(\tilde{\beta}^k)^T\bx_i = (\ln(k)/\mathrm{DegSEP}^*)y_i\bar\beta^T\bx_i \geq (\ln(k)/\mathrm{DegSEP}^*)\rho(\bar\beta) = \ln(k) \ .
\end{equation*}
Furthermore, since $\ell(t) \le \exp(-t)$, it holds that:
\begin{equation}\label{one_over_k}
\ell(y_i(\tilde{\beta}^k)^T\bx_i) \leq \ell(\ln(k)) \leq \exp(-\ln(k)) = 1/k \ .
\end{equation}
Clearly, the conditions for Proposition \ref{telgarsky_prop} are satisfied by $\ell_2$ steepest descent under the assumptions of Lemma \ref{telgarsky_norm_bound_lemma} (wherein $|S_k| = n$ for all $k$) as well as SGD under the assumptions of Lemma \ref{telgarsky_norm_bound_lemma_sgd} (wherein $|S_k| = 1$ for all $k$). (Note that in this proof $\alpha_j$ refers to the step-size with respect to the unnormalized version of $\ell_2$ steepest descent.)
Therefore, in both cases we may apply Proposition \ref{telgarsky_prop} using $\beta = \tilde\beta^k$ to yield:
\begin{equation*}
\|\beta^k - \tilde\beta^k\|_2^2 ~\leq~ \|\tilde\beta^k\|_2^2 + 2\sum_{j = 0}^{k - 1}\frac{\alpha_j}{|S_j|}\sum_{i \in S_j}\ell(y_i(\tilde{\beta}^k)^T\bx_i) ~\leq~ \frac{\ln(k)^2}{(\mathrm{DegSEP}^*)^2} + \frac{4}{\|\bX\|_{2,\infty}^2} \ ,
\end{equation*}
where the final inequality uses \eqref{one_over_k} as well as $\alpha_j \leq \frac{2}{\|\bX\|_{2,\infty}^2}$.
Therefore:
\begin{equation*}
\|\beta^k\| \leq \|\tilde\beta^k\|_2 + \|\beta^k - \tilde\beta^k\|_2 \leq \frac{\ln(k)}{\mathrm{DegSEP}^*} + \sqrt{\frac{\ln(k)^2}{(\mathrm{DegSEP}^*)^2} + \frac{4}{\|\bX\|_{2,\infty}^2}} \leq \frac{2\ln(k)}{\mathrm{DegSEP}^*} + \frac{2}{\|\bX\|_{2,\infty}} \ . 
\end{equation*}\ \qed

\subsection{Proof of Theorem \ref{sgdcomplexity}}\label{unitedflight}
To prove {\em (i)}, let $x \in \bbR^p$ be fixed and notice that for each $i \geq 0$ it holds that:
\begin{align*}
\|x^{i+1} - x\|_2^2 = \|x^{i} - \bar\alpha\tilde\nabla f(x^i) - x\|_2^2 = \|x^{i} - x\|_2^2 - 2\bar\alpha\tilde\nabla f(x^i)^T(x^i - x) + \bar\alpha^2\|\tilde\nabla f(x^i)\|_2^2
\end{align*} 
Rearranging terms, summing over $i \in \{0, \ldots, k\}$, and dividing by $2\bar\alpha(k+1)$ yields:
\begin{align}\label{classic_sequence_bound}
\frac{1}{k+1}\sum_{i = 0}^k \tilde\nabla f(x^i)^T(x^i - x) ~&=~ \frac{\|x^0 - x\|_2^2}{2\bar\alpha(k+1)} ~-~ \frac{\|x^{k+1} - x\|_2^2}{2\bar\alpha(k+1)} ~+~ \frac{\bar\alpha}{2(k+1)}\sum_{i = 0}^k \|\tilde\nabla f(x^i)\|_2^2 \nonumber \\
~&\leq~ \frac{\|x^0 - x\|_2^2}{2\bar\alpha(k+1)} ~+~ \frac{\bar\alpha}{2(k+1)}\sum_{i = 0}^k \|\tilde\nabla f(x^i)\|_2^2 \ .
\end{align}
Now, by the law of iterated expectations, for each $i \in \{0, \ldots, k\}$ it holds that
\begin{align*}
\bbE\left[\tilde\nabla f(x^i)^T(x^i - x)\right] &= \bbE\left[\bbE[\tilde\nabla f(x^i)^T(x^i - x) ~|~ x_i]\right] \\
&= \bbE\left[\bbE[\tilde\nabla f(x^i) ~|~ x_i ]^T(x^i - x)\right] \\
&= \bbE\left[\nabla f(x^i)^T(x^i - x)\right] \ ,
\end{align*}
where the last equality follows from the definition of the stochastic gradient, i.e., $\bbE[\tilde\nabla f(x^i) ~|~ x_i] = \nabla f(x^i)$. After taking the expectation of both sides of \eqref{classic_sequence_bound} and combining with the above we obtain
\begin{align*}
\frac{1}{k+1}\sum_{i = 0}^k \bbE\left[\nabla f(x^i)^T(x^i - x)\right] ~&\leq~ \frac{\|x^0 - x\|_2^2}{2\bar\alpha(k+1)} ~+~ \frac{\bar\alpha}{2(k+1)}\sum_{i = 0}^k \bbE[\|\tilde\nabla f(x^i)\|_2^2] \\
~&\leq~ \frac{\|x^0 - x\|_2^2}{2\bar\alpha(k+1)} ~+~ \frac{\bar\alpha M^2}{2} \ ,
\end{align*}
where the second inequality follows from \eqref{bounded_variance}. The gradient equality (which holds for each realization of $x^i$) states that $f(x^i) - f(x) \leq \nabla f(x^i)^T(x^i - x)$, and averaging the expectation of these inequalities over $i \in \{0, \ldots, k\}$ yields
\begin{equation}\label{final_exp_ineq}
\frac{1}{k+1}\sum_{i = 0}^k \bbE\left[f(x^i)\right] ~-~ f(x) ~\leq~ \frac{1}{k+1}\sum_{i = 0}^k \bbE\left[\nabla f(x^i)^T(x^i - x)\right] \ .
\end{equation}
Finally, in the case of Option A, Jensen's inequality implies {\em (i)} and, in the case of Option B, another iterated expectations argument implies that $\bbE[f(\hat{x}^k)] = \frac{1}{k+1}\sum_{i = 0}^k \bbE\left[f(x^i)\right]$ from which {\em (i)} directly follows.

Item {\em (ii)} follows directly from the format of the updates in Step (2.) of Algorithm \ref{sgd} as well as the triangle inequality. To prove {\em (iii)}, we use the smoothness of the objective function. In particular, applying \eqref{smooth-ineq} to the iterates of Algorithm \ref{sgd} yields for each $i \in \{0, \ldots, k\}$:
\begin{align}\label{smooth_bound_sgd}
f(x^{i+1}) ~&\leq~ f(x^{i}) + \nabla f(x^{i})^T(x^{i+1} - x^{i}) + \tfrac{L}{2}\|x^{i+1} - x^{i}\|^2_2 \nonumber \\
~&=~ f(x^{i}) - \bar\alpha\nabla f(x^{i})^T\tilde\nabla f(x^{i}) + \tfrac{L\bar\alpha^2}{2}\|\tilde\nabla f(x^{i})\|^2_2 \ . 
\end{align}
Notice that the law of iterated expectations as well as the definition of the stochastic gradient yields:
\begin{align*}
\bbE\left[\nabla f(x^{i})^T\tilde\nabla f(x^{i})\right] &= \bbE\left[\bbE[\nabla f(x^{i})^T\tilde\nabla f(x^{i}) ~|~ x^{i}]\right] \\
&= \bbE\left[\nabla f(x^{i})^T\bbE[\tilde\nabla f(x^{i}) ~|~ x_i]\right] \\
&= \bbE\left[\nabla f(x^{i})^T\nabla f(x^{i})\right] \\
&= \bbE\left[\|\nabla f(x^{i})\|_2^2\right] \ .
\end{align*}
Therefore, taking the expectation of both sides of \eqref{smooth_bound_sgd} and using \eqref{bounded_variance} yields:
\begin{equation*}
\bbE[f(x^{i+1})] ~\leq~ \bbE[f(x^{i})] - \bar\alpha \cdot \bbE\left[\|\nabla f(x^{i})\|_2^2\right] + \frac{L\bar\alpha^2 M^2}{2} \ .
\end{equation*}
Rearranging terms and summing over $i \in \{0, \ldots, k\}$ yields:
\begin{equation*}
\bar\alpha\sum_{i = 0}^k \bbE\left[\|\nabla f(x^{i})\|_2^2\right] ~\leq~ f(x^{0}) - \bbE[f(x^{k+1})] + \frac{L\bar\alpha^2 M^2 (k+1)}{2} \ .
\end{equation*}
Then using $f^\ast \leq \bbE[f(x^{k+1})]$ and dividing by $\bar\alpha(k+1)$ yields:
\begin{equation*}
\frac{1}{k+1}\sum_{i = 0}^k \bbE\left[\|\nabla f(x^{i})\|_2^2\right] ~\leq~ \frac{f(x^{0}) - f^\ast}{\bar\alpha(k+1)} + \frac{\bar\alpha L M^2}{2} \ .
\end{equation*}
Finally, since we are in the case of Option B, another iterated expectations argument implies that $\bbE\left[\|\nabla f(\hat{x}^k)\|_2^2\right] = \frac{1}{k+1}\sum_{i = 0}^k \bbE\left[\|\nabla f(x^{i})\|_2^2\right]$ from which {\em (iii)} directly follows. \ \qed

\subsection{Proof of Proposition \ref{LcalD}}\label{subway}
Recall that the scalar logistic loss function is denoted by $\ell(\cdot) : \bbR \to \bbR$, which is defined by $\ell(t) := \ln(1 + \exp(-t))$.
A simple calculation shows that $\ell^{\prime\prime}(t) = \frac{\exp(t)}{(\exp(t) + 1)^2} \leq \tfrac{1}{4}$ for all $t \in \bbR$.
As mentioned in Section \ref{sect:sgd-nonseparable}, it follows from item (2.) of Assumption \ref{Dassumption} that $L_\calD(\cdot)$ is continuous, convex, differentiable, and satisfies $\nabla L_\calD(\beta) = \bbE_{(\bx, y) \sim \calD}\left[\nabla_\beta \ell(y\beta^T\bx) \right] = \bbE_{(\bx, y) \sim \calD}\left[\ell^\prime(y\beta^T\bx) \cdot y\bx\right]$ for all $\beta \in \bbR^p$ (see, e.g., Section 7.2.4 of \cite{shapiro2009lectures}). Moreover, it can also be demonstrated (again by Section 7.2.4 of \cite{shapiro2009lectures}) that $L_\calD(\cdot)$ is twice differentiable. 
Letting $H(\beta)$ denote the Hessian matrix of $L_\calD(\beta)$ at $\beta$, then it holds that $H(\beta) = \bbE_{(\bx, y) \sim \calD}\left[\ell^{\prime\prime}(y\beta^T\bx)\bx\bx^T\right]$. 

Recall that a twice differentiable convex function is $L$-smooth with respect to the $\ell_2$ norm on $\bbR^p$ if and only if $H(\beta) \preceq LI_p$ for all $\beta \in \bbR^p$, where $I_p$ denotes the $p \times p$ identity matrix. Now, for any $(\bx, y)$, since $\ell^{\prime\prime}(y\beta^T\bx) \leq \tfrac{1}{4}$, we have that $\ell^{\prime\prime}(y\beta^T\bx)\bx\bx^T \preceq \tfrac{1}{4}\bx\bx^T$. Therefore, it holds that
\begin{equation*}
H(\beta) = \bbE_{(\bx, y) \sim \calD}\left[\ell^{\prime\prime}(y\beta^T\bx)\bx\bx^T\right] \preceq \tfrac{1}{4}\bbE_{(\bx, y) \sim \calD}\left[\bx\bx^T\right] = \tfrac{1}{4}\Sigma \preceq \tfrac{1}{4}\lambda_{\max}(\Sigma) I_{p} \ ,
\end{equation*}
which demonstrates that $L_\calD(\cdot)$ is $\tfrac{1}{4}\lambda_{\max}(\Sigma)$-smooth.

 \ \qed

\subsection{Proof of Proposition \ref{logit-var}}\label{valencia_cafesoret}
Again, recalling the notation $\ell(t) := \ln(1 + \exp(-t))$, note that $\ell^\prime(t) = -\tfrac{1}{\exp(t) + 1} \in (-1, 0)$ for all $t \in \bbR$. Then the stochastic gradient is $\nabla_\beta \ell(y\beta^T\bx)$ where $(\bx, y) \sim \calD$ and it holds that
\begin{align*}
\bbE_{(\bx, y) \sim \calD}\left[\|\nabla_\beta \ell(y\beta^T\bx)\|_2^2\right] &= \bbE_{(\bx, y) \sim \calD}\left[\|\ell^\prime(y\beta^T\bx) \cdot y\bx\|_2^2\right] \\
&= \bbE_{(\bx, y) \sim \calD}\left[|\ell^\prime(y\beta^T\bx)|\cdot\|\bx\|_2^2\right] \\
&\leq \bbE_{(\bx, y) \sim \calD}\left[\|\bx\|_2^2\right] \\
&= \tr(\Sigma) \ .
\end{align*} \ \qed

\end{document}